\documentclass[12pt]{article}

\usepackage{amssymb}
\usepackage{euscript}
\usepackage{mathrsfs} 
\usepackage[dvips]{graphicx}
\usepackage{flafter}
\usepackage[latin1]{inputenc}
\usepackage{amsmath}
\usepackage{amsthm}
\usepackage{amsfonts}
\usepackage{hyperref} 
\usepackage{amsxtra} 
\usepackage{verbatim} 
\usepackage{color} 
\usepackage[font={footnotesize}]{caption}  

\usepackage{pgf,tikz}
\usetikzlibrary{arrows}
\usepackage{caption}

\def\thesection{\arabic{section}}
\renewcommand{\theequation}{\thesection.\arabic{equation}}

\newtheorem{theorem}{Theorem}[section]
\newtheorem{lemma}[theorem]{Lemma}
\newtheorem{proposition}[theorem]{Proposition}

\newtheorem{question}[theorem]{Question}

\newtheorem{fact}[theorem]{Fact}

\theoremstyle{definition}   
\newtheorem{remark}[theorem]{Remark}

\evensidemargin\oddsidemargin

\newcommand{\eqnsection}{
\renewcommand{\theequation}{\thesection.\arabic{equation}}
    \makeatletter
    \csname  @addtoreset\endcsname{equation}{section}
    \makeatother}
\eqnsection

\def\P{{\bf P}}
\def\E{{\bf E}}

\def\ee{\mathrm{e}}

\def\TT{{\mathbb T}}

\begin{document}

 \baselineskip=18pt

\title{\bf Infinite  differentiability of the free energy  for a  Derrida-Retaux system}
\author{Xinxing Chen \\
\\
{\footnotesize School of Mathematical Sciences, Shanghai Jiaotong University, Shanghai 200240, China}\\
 {\footnotesize chenxinx@sjtu.edu.cn}
 } \maketitle

\bigskip
\bigskip
\bigskip

{\leftskip=1.5truecm \rightskip=1.5truecm \baselineskip=15pt \small

\noindent{\slshape\bfseries Abstract.} We consider a recursive system which was introduced by Derrida and Retaux (J. Stat. Phys. ${\bf 156}$ (2014) 268-290) as a toy model to study the depinning transition in presence of disorder. Derrida and Retaux predicted the free energy $F_\infty(p)$ of the system exhibit   quite an unusual physical phenomenon which is an  infinite order phase transition.
Hu and Shi (J. Stat. Phys. ${\bf 172}$ (2018) 718-741) studied  a special situation and obtained other behavior of the free energy, while insisted on $p=p_c$ being an essential singularity.
Recently,
 Chen, Dagard, Derrida, Hu, Lifshits and Shi (Ann. Probab. ${\bf 49}$ (2021) 637-670) confirmed the Derrida-Retaux conjecture under suitable integrability condition. However, in the mathematical review, it is still unknown whether the free  energy is infinitely differentiable  at the critical point.
 So that, we continue to  study the infinite differentiability of the free energy in this paper.

\bigskip


\bigskip


} 

\bigskip
\bigskip

\section{Introduction}

Fix  an integer $m\ge 2$. Let $X_0^*$ be a random variable which  takes    value of $\{1,2,3,\cdots\}$.
Let $p\in [0,1]$, and let $X_0$ be a nonnegative integer-valued random variable which satisfies $P_{X_0}=(1-p)\delta_0+pP_{X_0^*}$, i.e.,
\begin{equation}\label{e:lawX}
\P(X_0=0)=1-p,\quad \P(X_0=k)=p\P(X_0^*=k) \quad {\rm for~~ each~~} k\ge 1.
\end{equation}
For $n\ge 0$, we  recursively define
\begin{equation}\label{e:iterXXX}
X_{n+1}:=(X_{n,1}+X_{n,2}+\cdots+X_{n,m}-1)^+,
\end{equation}
where $X_{n,i}, i\ge 1$ are independent copies of $X_n$ and  $a^+=\max\{a,0\}$ for  $a\in {\mathbf{R}}$. To avoid trivialities, we may assume further $c_1:=\P(X_0^*\ge 2)> 0$ for otherwise $X_n\le 1, ~ n\ge 0$ almost surely when $m=2$.

 System $\{X_n, n\ge 0\}$ satisfied \eqref{e:iterXXX} is generally called the Derrida-Retaux system. Derrida-Retaux \cite{derrida-retaux} used it as a toy model to study the depinning transition,  Collet, Eckmann,   Glaser    and Martin  \cite{collet-eckmann-glaser-martin} as a
spin-glass model, Li and Rogers  \cite{li-rogers} as a hierarchical model. The Derrida-Retaux system is also   a max-type recursive distribution equation, see Aldous and Bandyopadhyay \cite{aldous-bandyopadhyay}. Moreover, it is  closely related to the parking model on an infinite regular tree, see Aldous, Contat, Curien and H$\acute{e}$nard \cite{aldous-Contat-curien-Henard}. There are difference continuous versions of Derrida-Retaux system, see Hu, Mallein and Pain \cite{HMP}, Chen, Dagard, Derrida and Shi \cite{4authors}. For more references and  conjectures about the Derrida-Retaux system, one can  see \cite{bz_DRsurvey}, \cite{xyz_sustainability} and \cite{xyz_dual}.

Collet, Eckmann,   Glaser    and Martin  \cite{collet-eckmann-glaser-martin} showed that there exists a phase transition for system $\{X_n, n\ge 0\}$.\\

\noindent {\bf{Theorem A}  }(Collet et al. \cite{collet-eckmann-glaser-martin})

 (1) If $(m-1)\E(X_0m^{X_0})\le \E(m^{X_0})<\infty$, then $\lim\limits_{n\rightarrow\infty}\E(X_n)=0$;

 (2) If $(m-1)\E(X_0m^{X_0})>\E(m^{X_0})$ or $\E(m^{X_0})=\infty$, then $\lim\limits_{n\rightarrow\infty}\E(X_n)=\infty$.\\

\noindent Indeed, as known in   \cite{collet-eckmann-glaser-martin} (also see \cite{bmxyz_questions}),
   $\E(m^{X_n})-(m-1)\E(X_n m^{X_n})$ keeps the same symbol  as    $\E(m^{X_0})-(m-1)\E(X_0 m^{X_0})$.
Accordingly, system $\{X_n,n\ge 0\}$ is said to be subcritical if  $ (m-1)\E(X_0m^{X_0})<\E(m^{X_0})  $, critical if $(m-1)\E(X_0m^{X_0})= \E(m^{X_0})<\infty$ and supercritical if $(m-1)\E(X_0m^{X_0})>\E(m^{X_0})$ or $\E(m^{X_0})=\infty$.
 Let
 $$p_c:=\frac{1}{1+\E(((m-1)X_0^*-1)m^{X_0^*})}\in [0,1).$$
The value $p=p_c$ is just  the unique solution    satisfying $\E(m^{X_0})=(m-1)\E(X_0 m^{X_0})$, provided $\E(X_0^* m^{X_0^*})<\infty$. So that, we can rewrite  Theorem A   as: Under $\E(X_0^*m^{X_0^*})<\infty$,
  $$
  \lim_{n\rightarrow\infty}\E(X_n)=0 \quad {\rm ~for~} \quad p\le p_c \quad\quad  {\rm  and } \quad\quad   \lim_{n\rightarrow\infty}\E(X_n)=\infty \quad  {\rm ~for~} \quad p> p_c.
  $$

It is  important to study  the quantity
 $$
 F_\infty(p):=\lim\limits_{n\rightarrow\infty}\downarrow\frac{\E(X_n)}{m^n}=\lim\limits_{n\rightarrow\infty}\uparrow \frac{\E(X_n)-\frac{1}{m-1}}{m^n},
 $$
which is called the free energy of system $\{X_n,n\ge 0\}$, see Derrida and Retaux \cite{derrida-retaux}.
By Theorem A, under $\E(X_0^*m^{X_0^*})<\infty$     the free energy also  has a phase transition:
$$
F_\infty(p)=0\quad {\rm   for} \quad p\le p_c \quad\quad {\rm ~and~~}\quad F_\infty(p)>0\quad {\rm  for } \quad p>p_c.
$$
Derrida and Retaux \cite{derrida-retaux} predicted system $\{X_n, n\ge 0\}$ exhibit an infinite order phase transition which is a Berezinskii-Kosterlitz-Thouless type, and gave a  famous conjecture which says that under suitable integrability condition on $X_0^*$, in the nearly supercritical regime
 $$
F_\infty(p)=\exp\{-\frac{C+o(1)}{(p-p_c)^{1/2}}\}, \quad p\downarrow p_c.
$$
Infinite order phase transition is quite an unusual physical phenomenon. Similar phenomenon were also shown in vertex-reinforced jump process on a regular tree \cite{peter-remy},  classical spin system on a lattice with a long range inhomogeneous coupling \cite{costin-costin-grunfeld} and
explosive percolation with a particular initial power-law distribution \cite{costa-dorogovtsev-goltsev-mendes}.

Later,  Hu and Shi \cite{yz_bnyz} considered a special $X_0^*$, which satisfies $\P(X_0^*=k)\sim c m^{-k} k^{-\alpha}$, $k\rightarrow\infty$, for some constants $c<\infty$ and   $\alpha\in {\mathbf{R}}$. They proved that if $\alpha<2$ then
\begin{equation}\label{e:expoDR1}
F_\infty(p)=\exp\{-\frac{1}{(p-p_c)^{\nu(\alpha)+o(1)}}\}, \quad p\downarrow p_c,
\end{equation}
where $ \nu(\alpha)=\frac{1}{2-\alpha}$ for $\alpha<2$ (note that $\E(X_0^*m^{X_0^*})=\infty$ in such situation).  Recently in \cite{CDDHLS2021},   we  gave a partial answer to the Derrida-Retaux conjecture by showing that if $\E((X_0^*)^3m^{X_0^*})<\infty$ then
\begin{equation}\label{e:expoDR2}
F_\infty(p)=\exp\{-\frac{1}{(p-p_c)^{1/2+o(1)}}\}, \quad p\downarrow p_c.
\end{equation}
 We also considered the same  situation  $\P(X_0^*=k)\sim c m^{-k} k^{-\alpha}$, $k\rightarrow\infty$,  and proved that \eqref{e:expoDR1} still  holds true  for $\alpha>2$, with $\nu(\alpha)=\begin{cases} \frac{1}{\alpha-2}, & 2<\alpha\le 4,\\ \frac{1}{2}, &\alpha>4.\end{cases}$ While when $\P(X_0^*=k)\sim c m^{-k} k^{-2}$, $k\rightarrow\infty$, the free energy was shown in a different behavior:
\begin{equation}\label{e:expoDR3}
F_\infty(p)=\exp(-e^{(c'+o(1))/p}),\quad p\downarrow p_c=0,
\end{equation}
where $c'=\frac{1}{(m-1)c}$.

However, in the mathematical review, it is still
  a problem whether  the free energy  $F_\infty(p)$ is infinitely differentiable at $p=p_c$.  Even in the situation $\P(X_0^*=k)\sim c m^{-k} k^{-\alpha}$, $k\rightarrow\infty$,   infinite differentiability  can not be derived directly from \eqref{e:expoDR1}, \eqref{e:expoDR2} or \eqref{e:expoDR3} owing to  the small term $o(1)$. We need some carefulness. As it is known, see  Russo \cite{russo}, the percolation probability on ${\mathbf{Z}}^2$ is  infinitely differentiable
     except  at $p=p_c$ at most.

Based on the above factors,  we  study the infinite differentiability  of $F_\infty(p)$ in this paper. Write $F_\infty^{(k)}(p)=\frac{d^k}{dp^k}F_\infty(p), ~p\in [0,1]$ for the $k$-th  derivative of $F_\infty(p)$, where $F_\infty^{(k)}(0)$ stands for the right derivative  at $p=0$, while $F_\infty^{(k)}(1)$ for the left derivative at $p=1$.

 By \eqref{e:iterXXX}, we have  $\E(X_{n+1})=m\E(X_n)-1+\P(X_n=0)^m$ which implies $\E(X_n)=m^n\E(X_0)-\frac{m^n-1}{m-1}+m^n \sum_{i=0}^{n-1}\frac{\P(X_i=0)^m}{m^{i+1}}$. So,
$$
F_\infty(p)=\lim_{n\rightarrow\infty}\frac{\E(X_n)}{m^n}=\E(X_0)-\frac{1}{m-1}+ \sum_{n=0}^{\infty}\frac{\P(X_n=0)^m}{m^{n+1}}.
$$
  Our main result is the following theorem.
\begin{theorem}\label{t:main}
Assume   $\E( s^{X_0^*})<\infty$ for  $|s|<m$.  Then  $F_\infty\in C^\infty[0,1]$ and for  $k\ge 0$,
$$
F_\infty^{(k)}(p)=\frac{d^k}{dp^k}(\E(X_0)-\frac{1}{m-1}  )+ \sum_{n=0}^{\infty}\frac{1}{m^{n+1}}\frac{d^k}{dp^k}\P(X_n=0)^m,
$$
where the summation converges uniformly in $p\in [0,1]$.
\end{theorem}

\begin{remark}\label{remark0} Let $\alpha\in \mathbf{R}$ and $c>0$. Consider the example $X_0^*$ which satisfies
 $\P(X_0^*=k)\sim c m^{-k} k^{-\alpha}$ as $k\rightarrow\infty$.  We have  $\E( s^{X_0^*})<\infty$ for  $|s|<m$ always,  and so $F_\infty\in C^\infty[0,1]$.
\end{remark}

\begin{remark}\label{remark2}  Assume   $\E( s^{X_0^*})<\infty$ for  $|s|<m$.  As the summation in Theorem \ref{t:main} converges uniformly, we have
$$
F_\infty^{(k)}(p)=\lim_{n\rightarrow\infty}\frac{d^k}{dp^k}\frac{\E(X_n)}{m^n}, \quad k\ge 0, \quad p\in [0,1].
$$
\end{remark}

\begin{remark}\label{remark1} Assume  $\E(X_0^*m^{X_0^*})<\infty$. Then $p_c\in (0,1)$.
 By Theorem \ref{t:main}, $F_\infty$ is infinitely differentiable at $p=p_c$. By Theorem A,  $F_\infty(p)=0$ for  $ p\in [0,  p_c)$.   Hence
$$
F_\infty^{(k)}(p_c)=\lim_{p\uparrow  p_c}F_\infty^{(k)}(p)=0,
\quad k\ge 0.
$$
\end{remark}

\bigskip

The rest of the paper is organized as follows. In Section 2, by using the Hoeffding's inequality and  admitting  Propositions \ref{t:pX1} and \ref{t:pX2}, we give the proof of  Theorem  \ref{t:main}. Where,
 Proposition \ref{t:pX1} gives upper bounds for $|\frac{d^k}{dp^k} \P(X_n=0)|$ in terms of $\P(X_n=0), \, \P(\sum\limits_{i=1}^{m^n-k} X_{0,i}\le m^n)$ and  $\prod\limits_{i=0}^{M-1}\E(~m^{(1-\delta)(X_i\wedge (M-i))}~)^{(m-1)k}$;     Proposition \ref{t:pX2} shows some inequalities related to  $\P(X_n=0)$ and  $\prod\limits_{i=0}^{M-1}\E(~m^{(1-\delta)(X_i\wedge (M-i))}~)^{m-1}$. In Sections 3 and 4,
we  prove   Propositions \ref{t:pX1} and \ref{t:pX2}   respectively. Some further remark and question is presented in Section 5.

 Notation. We will use $c_i>0, 1\le i\le 7$ and  $n_1, n_2\in \mathbf{Z}^+$    to stand for some constants which are independent of $p$ and $n$.

\section{Proof of Theorem \ref{t:main}}

Let us begin with the famous Hoeffding's inequality.

\begin{lemma}\label{HIE} (Hoeffding's inequality) Let $a>0$.   Let $W_n, n\ge 1$ be a sequence of independent random variables
with $0\le W_n\le a$.  Then
$$
\P(\sum_{i=1}^n W_i\le \sum_{i=1}^n \E(W_i)-t)\le e^{-\frac{2t^2}{na^2}}, \quad n\ge 1, \quad  t>0.
$$
\end{lemma}

By \eqref{e:iterXXX},  $\sum_{i=1}^m X_{n,i}-1 \le X_{n+1}\le \sum_{i=1}^m X_{n,i}$ for all $n\ge 0$. Hence in the meaning of
stochastic dominance, we have $\sum_{i=1}^{m^n} X_{0,i}-m^n \overset{st}{\le} X_n \overset{st}{\le} \sum_{i=1}^{m^n} X_{0,i}$,
where $X_{0,i}, i\ge 1$ are independent  copies of $X_0$. So that we need estimate the summation of these $X_{0,i}$. With the help of  Hoeffding's inequality, we have the following result.
\begin{lemma}\label{l:sum_m^n-k} Assume that  $c_1:=\P(X_0^*\ge 2)>0$.  Then for any $p\in [1-\frac{c_1}{4},1]$, $k\ge 1$ and $n\ge \lfloor \frac{\log(k)+\log(\frac{4+c_1}{c_1})}{\log m}\rfloor+1$, we have
$$
\P(\sum_{i=1}^{m^n-k} X_{0,i}\le m^n)\le e^{-\frac{c_1^2}{32} (m^n-k)}.
$$
\end{lemma}
 \proof Fix $p\in [1-\frac{c_1}{4},1]$, \, $k\ge 1$  and $n\ge \lfloor \frac{\log(k)+\log(\frac{4+c_1}{c_1})}{\log m}\rfloor+1$.
Then $(1+\frac{c_1}{4})(m^n-k)\ge m^n$.
   Since $c_1\in (0,1]$, \,  $p\in [1-\frac{c_1}{4},1]$,  \, $P_{X_0}=(1-p)\delta_0+pP_{X_0^*}$ and $X_0^*\in \{1,2,\cdots\}$,
$$
\E( X_0\wedge 2 )=p\E(X_0^*\wedge 2)= p(1+\P(X_0^*\ge 2))\ge (1-\frac{c_1}{4}) (1+c_1) \ge 1+\frac{c_1}{2}.
$$
Hence
$$
\sum_{i=1}^{m^n-k} \E(X_{0,i}\wedge 2)-\frac{c_1}{4}(m^n-k)= [\E( X_0\wedge 2 )- \frac{c_1}{4}](m^n-k) \ge m^n.
$$
Using Lemma \ref{HIE} we obtain immediately
\begin{align*}
\P(\sum_{i=1}^{m^n-k} X_{0,i}\le m^n)\le& \P\left(\sum_{i=1}^{m^n-k} (X_{0,i}\wedge 2)\le \sum_{i=1}^{m^n-k} \E(X_{0,i}\wedge 2)-\frac{c_1}{4}(m^n-k)\right)\\
\le&  e^{-\frac{c_1^2}{32}(m^n-k)}.
\end{align*}
\qed\\

Next, we give two propositions which will be useful in the proof of Theorem \ref{t:main}.

\begin{proposition}\label{t:pX1} Fix $\delta\in (0,\frac{1}{2})$.  Assume  $c_2:=\E(m^{(1-\delta)X_0^*})<\infty$. Let $p\in (0,1)$,  $k\ge 1$, $n\ge 1$  and $M=\lfloor n-\delta^{-1}\ln n\rfloor$. Then the following three statements hold true:

(1) $\Big|\frac{d^k}{dp^k} \P(X_n=0)\Big|\le 2^k \, k! \, m^{kn} \, \frac{\P(X_n=0)}{(1-p)^k}$;

(2) $\Big|\frac{d^k}{dp^k} \P(X_n=0)\Big|\le 2^k \, k! \, m^{kn} \,  \P(\sum_{i=1}^{m^n-k} X_{0,i}\le m^n)$;

(3) There exists  constant $n_1=n_1(m, \, \delta, \, c_2, \, k)\in \mathbf{Z}^+$  such that for  $n\ge n_1$,
$$
\Big|\frac{d^k}{dp^k} \P(X_n=0)\Big|\le m^{3\delta kn} \, \prod_{i=0}^{M-1}\E(~m^{(1-\delta)(X_i\wedge (M-i))}~)^{(m-1)k}.
$$
\end{proposition}

\begin{proposition}\label{t:pX2} Fix $\delta\in (0,\frac{1}{16m})$.   Assume that $c_1:=\P(X_0^*\ge 2)>0$ and  $c_2:=\E(m^{(1-\delta)X_0^*})<\infty$.  Then  there exist constants $c_i=c_i(m, \, \delta, \,  c_1,  \, c_2)>0, \, 3\le i\le 5$, such that  for  $p\in (0,1)$ and  $n>M\ge 1$,
\begin{equation}\label{e:px2}
\prod_{i=0}^{M-1}\E(m^{(1-\delta)(X_i\wedge (M-i))})^{m-1}\le   c_3m^{2\delta M} \quad  {\rm or} \quad \P(X_{n}=0)\le  c_5e^{-c_4 m^{n-M}}.
\end{equation}
\end{proposition}

Let us admit   Propositions \ref{t:pX1} and \ref{t:pX2} hold true for the time being  whose proof will be postponed to Sections 3 and 4 respectively,
and we will use them to prove Theorem \ref{t:main}. \\

{\it Proof of Theorem \ref{t:main}.} By $X_{n+1}=(X_{n,1}+\cdots+X_{n,m}-1)^+$ for $n\ge 0$, we have $\E(X_{n+1})=m\E(X_n)-1+\P(X_n=0)^m$ and so,
$$
\E(X_n)=m^n\E(X_0)-\frac{m^n-1}{m-1}+\sum_{i=0}^{n-1}m^{n-i-1}\P(X_i=0)^m.
$$
Hence
$$
F_\infty(p)=\lim_{n\rightarrow\infty}\frac{\E(X_n)}{m^n}=\E(X_0)-\frac{1}{m-1}+\sum_{n=0}^{\infty}\frac{\P(X_n=0)^m}{m^{n+1}}.
$$

Since $P_{X_0}=(1-p)\delta_0+pP_{X_0^*}$,  $\P(X_0=\ell)$ is a polynomial function of $p\in [0,1]$ for each $\ell\ge 0$. Since  $X_{n+1}=(X_{n,1}+\cdots+X_{n,m}-1)^+$ for all $n\ge 0$, we can iteratively get  that $\P(X_n=\ell)$ is a polynomial function of $p$ for each $n\ge 0$ and $  \ell\ge 0$, too. It  deduces that all $\P(X_n=\ell)$ are infinitely differentiable in  $p$; Especially they have  right derivatives at $p=0$ and left derivatives at $p=1$. So, we will have our main result
if   for each $k\ge 0$ there has
\begin{equation}\label{e:convergence}
\lim_{n_0\rightarrow\infty}\sup_{p\in (0,1)}\sum_{n=n_0}^\infty \frac{1}{m^{n+1}}\left|\frac{d^k}{dp^k}\P(X_n=0)^m\right|=0.
\end{equation}
Since $\frac{d^k}{dp^k}\P(X_n=0)^m=\sum\limits_{(k_1, ~k_2, \cdots, ~k_m): ~k_1 + k_2 + \cdots + k_m=k}\binom{k}{k_1, ~k_2, \cdots, ~k_m} \prod\limits_{i=1}^m \frac{d^{k_i}}{dp^{k_i}}\P(X_n=0)$, we have
$$
\left|\frac{d^k}{dp^k}\P(X_n=0)^m\right|\le m^k \max\left\{|\P(X_n=0)|, |\frac{d}{dp} \P(X_n=0)|, \cdots, |\frac{d^k}{dp^k} \P(X_n=0)|\right\}^m .
$$
So, to prove  \eqref{e:convergence}, it is suffice to prove that for each $k\ge 0$ uniformly in $p\in (0,1)$ as $n\rightarrow\infty$,
\begin{equation}\label{e:derivatives}
\frac{d^k}{dp^k}\P(X_n=0)=O(m^{\frac{n}{2m}}).
\end{equation}

   When $k=0$, it is trivial for \eqref{e:derivatives}  since $0\le \P(X_n=0)\le 1$.

   Fix $k\ge 1$ and $\delta\in (0,\frac{\ln m}{10 k m^2})$.  Let $n_1\in \mathbf{Z}^+$ and   $c_i>0, \, 3\le i\le 5$ be the constants in Propositions \ref{t:pX1} and \ref{t:pX2}.
   Let $n\ge \max\{  \lfloor \frac{\log(k)+\log(\frac{4+c_1}{c_1})}{\log m}\rfloor+1, n_1\}$.
   If $p\in [1-\frac{c_1}{4},1)$, then using (2) of Proposition \ref{t:pX1}  and Lemma \ref{l:sum_m^n-k},  we get
$$
   \Big|\frac{d^k}{dp^k} \P(X_n=0)\Big|\le 2^k k! m^{kn}  \P(\sum_{i=1}^{m^n-k} X_{0,i}\le m^n)\le 2^k k! m^{kn} e^{-\frac{c_1^2}{32} (m^n-k)}=O(m^{\frac{n}{2m}}).
$$
  Otherwise, let $p\in (0,1-\frac{c_1}{4})$.  By (1) and (3) of Proposition \ref{t:pX1},   we   have
 \begin{align*}
   \Big|\frac{d^k}{dp^k} \P(X_n=0)\Big|\le \min\{ 2^k k! m^{kn} \frac{\P(X_n=0)}{(\frac{c_1}{4})^k}, ~  m^{3\delta kn} \prod_{i=0}^{M-1}\E(m^{(1-\delta)X_i^{(M)}})^{(m-1)k}\},
   \end{align*}
   where $M=\lfloor n-\delta^{-1}\ln n\rfloor$ and $X_i^{(M)}=X_i\wedge (M-i)$. By Proposition \ref{t:pX2}, there has
  $$
\prod_{i=0}^{M-1}\E(m^{(1-\delta)X_i^{(M)}})^{m-1}\le   c_3m^{2\delta M} \quad\quad {\rm or} \quad\quad  \P(X_{n}=0)\le  c_5e^{-c_4 m^{n-M}}.
$$
Hence
\begin{align*}
 \Big|\frac{d^k}{dp^k} \P(X_n=0)\Big|\le& \max\{ 2^k k! m^{kn} \frac{ c_5e^{-c_4 m^{n-M}}      }{(\frac{c_1}{4})^k}, \,  m^{3\delta kn} \big(c_3m^{2\delta M}   \big)^{k}\}\\
 \le & \max\{8^k k! c_5c_1^{-k} m^{kn} e^{-c_4 m^{\delta^{-1} \ln n}}, \,  c_3^{k} m^{5\delta kn}\}.
\end{align*}
Since $\delta\in (0,\frac{\ln m}{10 km^2}), \, k\ge 1$ and $m\ge 2$, we have $\delta^{-1}  \ge \frac{40}{\ln m}$ and $5\delta k\le \frac{1}{2m}$. Therefore, uniformly in $p\in (0,1-\frac{c_1}{4})$ as $n\rightarrow\infty$,
$$
\Big|\frac{d^k}{dp^k} \P(X_n=0)\Big|\le \max\{8^k k! c_5c_1^{-k} m^{kn}e^{-c_4 n^{40}}, \, c_3^{k} m^{\frac{n}{2m}}\}=O(m^{\frac{n}{2m}}).
$$
Such we prove \eqref{e:derivatives} and finish the proof of the theorem.  \qed\\

\section{Proof of Proposition \ref{t:pX1}}
To obtain the bounds of $\frac{d^k}{dp^k}\P(X_n=0)$ in Proposition \ref{t:pX1}, it is convenient to use  a  hierarchical representation of system $\{X_n, n\ge 0\}$,
as in \cite{collet-eckmann-glaser-martin,derrida-retaux,CDDHLS2021}.

 Let  $\mathbb{T}$ be a (reversed) $m$-regular tree.
 For any  vertex $v\in \mathbb{T}$, denote by $|v|$  the generation of $v$. Let $\TT_n :=\{v\in \TT: |v|=n\}$ for  $n\ge 0$.  So that, the initial generation $\TT_0$ is just the set of  the leaves of $\mathbb{T}$. For  $v\in \TT\backslash \TT_0$,  let $v^{(1)}, \, v^{(2)},\cdots, \, v^{(m)}$ be the $m$ parents of $v$.

  For  $v\in \TT_0$,  let $X_0^*(v)$ be a random variable having the law as $X_0^*$,  $U(v)$  a binomial random variable with  $\P(U(v)=1)=p$.  Assume further all these $X_0^*(v), \, U(v), \, v\in \TT_0$ are independent. Define
\begin{equation}\label{e:X(v)}
X(v) :=X_0^*(v)U(v), \quad v\in \TT_0;
\end{equation}
and iteratively set
\begin{equation}\label{e:hirachical_rep}
X(v) :=(  X(v^{(1)})+\cdots+ X(v^{(m)})-1)^+, \quad v\in \mathbb{T}\setminus\mathbb{T}_0.
\end{equation}
By definition,   $X(v), v\in \TT_n$ are independent and identically distributed (i.i.d) having  the same law as $X_n$. We will use $X$ to stand for  $(X(v), v\in \mathbb{T})$.

 For $n\ge 0$, write $\ee_n$ for the  first  lexicographic vertex in the $n$-th generation of $\TT$. Since $X(\ee_n)$ has the same law as $X_n$, we will use the notation $X_n=X(\ee_n)$ if without making  any confusion.
 For  $u\in \mathbb{T}$ and $0\le n\le |u|$, set $\mathbb{T}^u :=\{v\in \mathbb{T}: v $ is  an   ancestor of $u\}\cup\{u\}$  and
  $\mathbb{T}^{u}_n := \mathbb{T}^u\cap \TT_n$. 
Then the value of  $X_n$ is determined by   $X(v), v\in \mathbb{T}_0^{\ee_n}$.
 Therefore, our question is changed into  calculating the derivatives of $\E(f(X))$, where $f:\mathbf{R}^{\TT} \rightarrow \mathbf{R}$ is some indicator function    satisfying that the value of $f(X)$ is determined by $X(v), v\in \mathbb{T}_0^{\ee_n}$.

To obtain  a simple form of our result, we still need some notation. Let $A$ be a finite subset of $\TT_0$. We define  $\Theta^A X:=(\Theta^A X(v), v\in \mathbb{T})$ which such that
\begin{align}
\label{e: hierarchical0*} \Theta^A X(v) :=&X(v)1_{\{v\not\in A\}}+X_0^*(v)1_{\{v\in A\}}, ~~v\in \TT_0;\\
\label{e: hierarchical}\Theta^A X(v) :=&(\Theta^A X(v^{(1)})+\cdots+\Theta^A X(v^{(m)})-1)^+, ~~v\in \TT\backslash \TT_0.
\end{align}
Then the value of $\Theta^A X(\ee_n)$ is determined by $\Theta^A X(v), \, v\in \TT_0^{\ee_n}$ for each $n$.
Since $X(u)\le X_0^*(u)$ for $u\in \TT_0$, we have $\sup_{B:B\subset A}\Theta^B X(v)=\Theta^A X(v)$  for any $v\in \TT_0$. Iteratively using \eqref{e: hierarchical},  we have $\sup_{B:B\subset A}\Theta^B X(v)=\Theta^A X(v)$ for  $v\in \TT \backslash \TT_0$, too.
For   any function $f$ on  $\mathbf{R}^{\mathbb{T}}$, we define
$$
\nabla^A f(X) := \sum_{B: B\subseteq A} (-1)^{|A|-|B|}  f(\Theta^B X),
$$
where $|A|$ is the cardinality of $A$. By definition, 
%
  for any constants $ a,b\in {\mathbf{R}}$ and functions $f$ and $g$ on ${\mathbf{R}}^\TT$,
\begin{equation}\label{e:nabla+}
\nabla^A(af+bg)(X)=a\nabla^A f(X)+b\nabla^A g(X).
\end{equation}
%




\begin{lemma}\label{l:dev_on_tree}  Let $k\ge 1$, \, $n\ge 0$ and $x_v\in \mathbf{Z}^+$ for $v\in \TT_0^{\ee_n}$. Then
 \begin{align}\label{e:dev_on_tree_e}
\frac{d^k}{dp^k} \P(&X(v)=x_v, v\in \mathbb{T}_0^{\ee_n})
 = \frac{k!}{(1-p)^k} \sum_{A\subset  \TT_0^{\ee_n}:|A|=k} \E(1_{\{ X|_{A}=0\}}\nabla^A 1_{\{ X(v)=x_v, \, v\in \mathbb{T}_0^{\ee_n}  \}}   ),
\end{align}
where the meaning of $X|_A=0$ is  $X(v)=0$ for all $v\in A$.
\end{lemma}
\proof  Fix $x_v\in {\mathbf{Z}}^+, \, v\in \mathbb{T}_0^{\ee_n}$ and set ${\mathcal{D}} :=\{v\in \mathbb{T}_0^{\ee_n}: x_v>0\}$, $\alpha_{x,\mathcal{D}} :=\prod_{v\in {\mathcal{D}}} \P(X_0^*(v)=x_v)$.

  The left-hand side of \eqref{e:dev_on_tree_e} is easy to calculate. By \eqref{e:X(v)}, for  $v\in \mathbb{T}_0^{\ee_n}$,
 $$\P(X(v)=x_v)=\P( X_0^*(v)U(v)=x_v )=[\P(X_0^*(v)=x_v)p]^{1_{\{x_v>0\}}}(1-p)^{1_{\{x_v=0\}}}.$$
  So,
\begin{align*}\P(X(v)=x_v, v\in \mathbb{T}_0^{\ee_n})=&\prod_{v\in \mathbb{T}_0^{\ee_n}}\P(X(v)=x_v)=\alpha_{x,\mathcal{D}}~p^{ |{\mathcal{D}}|}~(1-p)^{m^n-|{\mathcal{D}}|} .
\end{align*}
 As a result,
\begin{align*}
  {\rm LHS}_{\eqref{e:dev_on_tree_e}}=& \alpha_{x,\mathcal{D}} \frac{d^k}{dp^k}  p^{ |{\mathcal{D}}|}  (1-p)^{m^n-|{\mathcal{D}}|} \\
 =&\alpha_{x,\mathcal{D}}\sum_{h=0}^{k}\binom{k}{h} (\frac{d^h}{dp^h}p^{|{\mathcal{D}}|})(\frac{d^{k-h}}{dp^{k-h}}(1-p)^{m^n-|{\mathcal{D}}|}) \\
 =&k!\alpha_{x,\mathcal{D}}\sum_{h=0}^{k}(-1)^{k-h}\binom{|{\mathcal{D}}|}{h}\binom{m^n-|{\mathcal{D}}|}{k-h}   p^{|{\mathcal{D}}|-h} (1-p)^{m^n-|{\mathcal{D}}|-k+h}.
\end{align*}

Next, we calculate the right-hand side  of \eqref{e:dev_on_tree_e}. Fix $A\subset  \TT_0^{\ee_n}$ with $|A|=k$. Let $B\subset A$.
By \eqref{e: hierarchical0*},  conditioned on  event $\{X|_A=0,~ 1_{\{ \Theta^B X(v)=x_v, v\in \mathbb{T}_0^{\ee_n}  \}}=1\}$,  for   $u\in A\cap \mathcal{D}$,
$$
1_{\{u\in B\}}\ge 1_{\{X(u)=0, ~\Theta^B X(u)>0\}}\ge 1_{\{X|_A=0, ~\Theta^B X(u)=x_u\}}=1;
$$
while  $x_u=\Theta^B X(u)=X_0^*(u)\ge 1$ for   $u\in B$. Hence there must has
 $B=A\cap {\mathcal{D}}$ conditioned on   $\{1_{\{X|_A=0\}} 1_{\{ \Theta^B X(v)=x_v, v\in \mathbb{T}_0^{\ee_n}  \}}=1\}$. So, writing $h=|A\cap {\mathcal{D}}|$, we have
 \begin{align*}
1_{\{ X|_{A}=0\}}\nabla^A 1_{\{ X(v)=x_v, v\in \mathbb{T}_0^{\ee_n}  \}}
=&1_{\{ X|_{A}=0\}} \sum_{B:B\subset A} (-1)^{|A|-|B|}  1_{\{\Theta^B  X(v)=x_v, v\in \mathbb{T}_0^{\ee_n}  \}}\\
=& 1_{\{ X|_{A}=0\}}~  (-1)^{k-h}  1_{\{ \Theta^{A\cap {\mathcal{D}}} X(v)=x_v, v\in \mathbb{T}_0^{\ee_n}  \}}\\
=&(-1)^{k-h} \prod_{v\in {\mathcal{D}}}1_{\{X_0^*(v)=x_v\}} \prod_{v\in {\mathcal{D}}\setminus A} 1_{\{U(v)=1\}} \prod_{v\in \TT_0^{\ee_n}\setminus({\mathcal{D}}\setminus A) } 1_{\{U(v)=0\}}.
\end{align*}
It follows immediately,
\begin{align*}
&\E(1_{\{ X|_{A}=0\}}\nabla^A 1_{\{ X(v)=x_v, v\in \mathbb{T}_0^{\ee_n}  \}})\\
=&(-1)^{k-h} \prod_{v\in {\mathcal{D}}}\P(X_0^*(v)=x_v) ~ \prod_{v\in {\mathcal{D}}\setminus A} \P(U(v)=1)~\prod_{v\in \TT_0^{\ee_n}\setminus({\mathcal{D}}\setminus A) } \P(U(v)=0)\\
=&(-1)^{k-h}\alpha_{x,\mathcal{D}}~p^{|{\mathcal{D}}|-h}(1-p)^{m^n-|{\mathcal{D}}|+h}.
\end{align*}
Therefore,
\begin{align*}
{\rm RHS}_{\eqref{e:dev_on_tree_e}}
=&\frac{k!}{(1-p)^k}\sum_{h=0}^k \sum_{A\subset  \TT_0^{\ee_n}:|A|=k, |A\cap {\mathcal{D}}|=h} \E(1_{\{ X|_{A}=0\}}\nabla^A 1_{\{ X(v)=x_v, \, v\in \mathbb{T}_0^{\ee_n}  \}}   )\\
=& \frac{k!}{(1-p)^k} \alpha_{x,\mathcal{D}}\sum_{h=0}^k \sum_{A\subset  \TT_0^{\ee_n}:|A|=k, |A\cap {\mathcal{D}}|=h} (-1)^{k-h}~p^{|{\mathcal{D}}|-h}(1-p)^{m^n-|{\mathcal{D}}|+h}\\
=&   k!  \alpha_{x,\mathcal{D}}\sum_{h=0}^k  \binom{|{\mathcal{D}}|}{h}\binom{m^n-|{\mathcal{D}}|}{k-h} (-1)^{k-h}~p^{|{\mathcal{D}}|-h}(1-p)^{m^n-|{\mathcal{D}}|-k+h}.
\end{align*}
Hence   ${\rm LHS}_{\eqref{e:dev_on_tree_e}}=  {\rm RHS}_{\eqref{e:dev_on_tree_e}}$ holds true, we complete the proof of the lemma.
 \qed\\

\begin{lemma}\label{l:Dev_F_n} For each $k\ge 1$ and $n\ge 0$, we have
$$
|\frac{d^k}{dp^k}\P(X_n=0)|\le  \frac{2^k k!}{(1-p)^k}\sum_{A\subset \mathbb{T}_0^{\ee_n}:~|A|=k}  \P(\nabla^{A} 1_{\{X_n=0\}}\not=0,~ X|_A=0).
$$
\end{lemma}

\proof Fix $k\ge 1$ and $n\ge 0$. By \eqref{e:hirachical_rep}, $X(v)\ge X(v^{(j)})-1$ for  $\TT\setminus \TT_0$ and $1\le j\le m$. Hence
$$X(\ee_n)\ge X(v)-n+|v|, \quad v\in \TT^{\ee_n}.$$
It implies that conditioned on $\{X_n=0\}$ we have $X(v)\le n$ for all $v\in \TT_0^{\ee_n}$. Hence
$$
1_{\{X_n=0\}}=\sum_{x_v\in\{0,1,\cdots,n\} {\rm ~for~} v\in \TT_0^{\ee_n}} 1_{\{X_n=0, ~X(v)=x_v, ~v\in \mathbb{T}_0^{\ee_n}\}}.
$$
 Since the value of $X_n$ is determined by $X(v), v\in \mathbb{T}_0^{\ee_n}$, the above equation can be rewrote as
\begin{equation}\label{e:XXXnnn}
 1_{\{X_n=0\}}=\sum_{x_v\in\{0,1,\cdots,n\} {\rm ~for~} v\in \TT_0^{\ee_n}} a_x 1_{\{X(v)=x_v, ~v\in \mathbb{T}_0^{\ee_n}\}},
 \end{equation}
where $a_x\in \{0,1\}$ is non-random and satisfies $1_{\{X_n=0, X(v)=x_v, v\in \TT_0^{\ee_n}\}}=a_x 1_{\{ X(v)=x_v, v\in \TT_0^{\ee_n}\}}$.  So,
\begin{equation}\label{e:finite_term}
\P(X_n=0)=\sum_{x_v\in\{0,1,\cdots,n\} {\rm ~for~} v\in \TT_0^{\ee_n}} a_x~ \P(X(v)=x_v, ~v\in \mathbb{T}_0^{\ee_n}) .
\end{equation}
Since  the sums on the right-hand side of \eqref{e:finite_term}   have only finite terms,  we have
$$
\frac{d^k}{dp^k}\P(X_n=0)= \sum_{x_v\in\{0,1,\cdots,n\} {\rm ~for~} v\in \TT_0^{\ee_n}} a_x~\frac{d^k}{dp^k }\P(X(v)=x_v, ~v\in \mathbb{T}_0^{\ee_n}).
$$
Using Lemma \ref{l:dev_on_tree} and the Fubini Theorem, we  get
\begin{align*}
\frac{d^k}{dp^k}\P(X_n=0)=& \sum_{x_v\in\{0,1,\cdots,n\} {\rm ~for~} v\in \TT_0^{\ee_n}} a_x \frac{k!}{(1-p)^k} \sum_{A\subset  \TT_0^{\ee_n}:|A|=k} \E(1_{\{ X|_{A}=0\}}\nabla^A 1_{\{ X(v)=x_v, v\in \mathbb{T}_0^{\ee_n}  \}}   )\\
=& \frac{k!}{(1-p)^k} \sum_{A\subset  \TT_0^{\ee_n}:|A|=k}   \E(1_{\{ X|_{A}=0\}} \sum_{x_v\in\{0,1,\cdots,n\} {\rm ~for~} v\in \TT_0^{\ee_n}} a_x \nabla^A 1_{\{ X(v)=x_v, v\in \mathbb{T}_0^{\ee_n}  \}}   ).
\end{align*}

On the other hand, by \eqref{e:XXXnnn} and \eqref{e:nabla+},
\begin{align*}
\sum_{x_v\in\{0,1,\cdots,n\} {\rm ~for~} v\in \TT_0^{\ee_n}} a_x \nabla^A 1_{\{ X(v)=x_v, v\in \mathbb{T}_0^{\ee_n}  \}} =&\nabla^A  \sum_{x_v\in\{0,1,\cdots,n\} {\rm ~for~} v\in \TT_0^{\ee_n}} a_x  1_{\{ X(v)=x_v, v\in \mathbb{T}_0^{\ee_n}  \}}\\
=&\nabla^A 1_{\{X_n=0\}}.
\end{align*}
Therefore,
$$
\frac{d^k}{dp^k}\P(X_n=0)=\frac{k!}{(1-p)^k} \sum_{A\subset  \TT_0^{\ee_n}:|A|=k}   \E(1_{\{ X|_{A}=0\}} \nabla^A 1_{\{X_n=0\}}).
$$

  Since    $|\nabla^A 1_{\{X_n=0\}}|= |\sum\limits_{B:B\subset A} (-1)^{|A|-|B|} 1_{\{ \Theta^B X(\ee_n)=0\}}|\le 2^k$ for   $A\subset \mathbb{T}_0^{\ee_n}$ with $|A|=k$, we draw out the conclusion of
 the  lemma     immediately.\qed\\

We have to estimate $\P(\nabla^{A} 1_{\{X_n=0\}}\not=0,~ X|_A=0)$, so that, we should understand what will happen when event $\{\nabla^{A} 1_{\{X_n=0\}}\not=0\}$ occurs.
 \begin{lemma}\label{l:AXs00} Let $n\ge 0, \, i\ge 0$ and $A \subset \mathbb{T}^{\ee_n}_0$ with $A\not=\emptyset$. Then conditioned on event    $\{\nabla^A 1_{\{X(\ee_n)=i\}}\not=0\}$,   the following three statements hold true:
\begin{align*}
(1)&~X(v)\le n+i-|v| {\rm ~for~ any~} v\in \mathbb{T}^{\ee_n} ;\\
(2)&~\Theta^A X(\ee_n)\ge (i\vee 1);\\
(3)&~{\rm There~  exist ~ some~integers~} x_1\ge 0, \cdots, x_m\ge 0 {\rm ~with~}(x_1+\cdots+x_m-1)^+=i {\rm ~such~that}\\
&\quad\quad \quad\quad \quad\quad \quad\quad \nabla^{A\cap\mathbb{T}_0^{\ee_n^{(j)}}} 1_{\{ X(\ee_{n}^{(j)})=x_j\}}\not=0 \quad {\rm for~each~}1\le j\le m.
\end{align*}
\end{lemma}
 \proof Fix $n\ge 0, \, i\ge 0$ and $A \subset \mathbb{T}^{\ee_n}_0$ with $A\not=\emptyset$.  Suppose  event     $\{\nabla^A 1_{\{X(\ee_n)=i\}}\not=0\}$ occur.

Since $\nabla^A 1_{\{X(\ee_n)=i\}}=\sum\limits_{B:B\subset A} (-1)^{|A|-|B|} 1_{\{ \Theta^B X(\ee_n)=i\}}$ and $\sum\limits_{B:B\subset A} (-1)^{|A|-|B|}=(1-1)^{|A|}=0$,
$$
0=\inf_{B:B\subset A}1_{\{ \Theta^B X(\ee_n)=i\}}<\sup_{B:B\subset A}1_{\{ \Theta^B X(\ee_n)=i\}}=1.
$$
Since  $\Theta^B X(\ee_n)\ge X(\ee_n)$ for any $B$,
$$
1_{\{X(\ee_n)\le i\}}=\sup_{B:B\subset A}1_{\{ \Theta^B X(\ee_n)\le i\}}\ge  \sup_{B:B\subset A}1_{\{ \Theta^B X(\ee_n)=i\}}.
$$
Hence $X(\ee_n)\le i$. Furthermore, by $X(\ee_n)\ge X(v)-n+|v|$ for $v\in \mathbb{T}^{\ee_n}$,  we obtain (1).

Since  $\sup_{B: B\subset A} \Theta^B X(\ee_n)= \Theta^A X(\ee_n)$, we have
$$
1_{\{i=0\}}1_{\{\Theta^A X(\ee_n)=0\}}\le \inf_{B:B\subset A}1_{\{ \Theta^B X(\ee_n)=i\}}=0.
$$
and
$$
1_{\{\Theta^A X(\ee_n)\ge i\}}=\sup_{B:B\subset A}1_{\{ \Theta^B X(\ee_n)\ge i\}}\ge  \sup_{B:B\subset A}1_{\{ \Theta^B X(\ee_n)=i\}}=1.
$$
So that we draw out the conclusion of   (2).

We are left to prove (3). For short, we write $[m]=\{1,2,\cdots,m\}$,  ~~$A_{j}= A\cap \mathbb{T}_0^{\ee_{n}^{(j)}}$ for $j\in [m]$ and
$\Xi_i=\{(x_1,\cdots,x_m): (x_1+\cdots+x_m-1)^+=i, \, x_j\ge 0 {\rm ~for~}   j\in [m]\}$. By
 \eqref{e:nabla+} and   $1_{\{X(\ee_{n})=i\}}=\sum\limits_{(x_1,\cdots,x_m)\in \Xi_i}  1_{\{X(\ee_n^{(j)})=x_j ~{\rm for} ~j\in [m]\}}$, we  have
 $$
 \nabla^A 1_{\{X(\ee_{n})=i\}}=\sum_{(x_1,\cdots,x_m)\in \Xi_i}\nabla^A   1_{\{X(\ee_{n}^{(j)})=x_j ~{\rm for} ~j\in [m]\}}.
 $$
  On the other hand,
   since $\Theta^{B_1\cup\cdots \cup B_m} X(\ee_n^{(j)})= \Theta^{B_j} X(\ee_n^{(j)})$ for any $B_1\subset A_1, \cdots, ~B_m\subset A_m$ and $j\in [m]$, we get that  for each $(x_1,\cdots, x_m)\in \Xi_i$,
 \begin{align*}
 \nabla^A   1_{\{X(\ee_{n}^{(j)})=x_j ~{\rm for}~j\in [m]\}}=&\sum_{B_1\subset A_1,\cdots,B_m\subset A_m} (-1)^{|A|-|B_1\cup \cdots \cup B_m|}  1_{\{\Theta^{B_1\cup\cdots \cup B_m}  X(\ee_{n}^{(j)})=x_j {\rm ~for~} j\in [m]\}}\\
 =&\sum_{B_1\subset A_1,\cdots,B_m\subset A_m}(-1)^{(|A_1|+\cdots+|A_m|)-(|B_1|+ \cdots+ B_m|)}  1_{\{\Theta^{B_j}  X(\ee_{n}^{(j)})=x_j {\rm ~for~} j\in [m]\}}\\
 =&\prod_{j=1}^m  \sum_{B_j\subset A_j } (-1)^{|A_j|-|B_j|}  1_{\{\Theta^{B_j}  X(\ee_{n}^{(j)})=x_j\}}\\
 =& \prod_{j=1}^m  \nabla^{A_j}   1_{\{X(\ee_{n}^{(j)})=x_j\}}.
 \end{align*}
 Therefore,
 $$
 \nabla^A 1_{\{X(\ee_{n})=i\}}=\sum_{(x_1,\cdots,x_m)\in \Xi_i}\prod_{j=1}^m  \nabla^{A_j}   1_{\{X(\ee_{n}^{(j)})=x_j\}}.
 $$
Hence we have (3). \qed\\

Let us introduce some  notation again. For $v\in \mathbb{T}_0$ and $n\ge 1$, let  $v_n$ denote  the unique descendant of $v$ in generation $n$.
 For $u\in \TT$ and $A\subset\TT_0^{u}$, set
  $$
\mathbb{O}_{u,A} :=\{v_i:  v\in A, \, 1\le i\le |u| \},
$$
$$
\mathbb{L}_{u,A} :=\{v^{(j)}:  v\in \mathbb{O}_{u,A}, \, 1\le j\le m   \}\setminus (\mathbb{O}_{u,A}\cup A).
$$
Then
 $A\cup \mathbb{O}_{u,A}$ is the set of vertices of the smallest subtree of $\mathbb{T}$ which contains $A\cup\{u\}$  if $|u|\ge 1$, while $\mathbb{O}_{u,A}=\emptyset$ if $|u|=0$. By observed,   all $\TT^v$, $v\in A\cup \mathbb{L}_{u,A}$ are disjoint subsets of $\TT^{u}$.
  Since  $v_{i}\in \{v_{i+1}^{(j)}: 1\le j\le m\}$ for any $v\in \TT_0$ and $i\ge 0$, we have
$$| \mathbb{L}_{u,A}\cap \TT_i|\le (m-1) |A|, \quad i\ge 0.$$
\begin{lemma}\label{l:AXs} Let $n\ge 0$. Then
 for each $i\ge 0$  and $A \subset \mathbb{T}^{\ee_n}_0$ with $A\not=\emptyset$,
 conditioned on   event    $\{\nabla^A 1_{\{X(\ee_n)=i\}}\not=0\}$,
\begin{equation}\label{e:AXs}
\Theta^A X(\ee_n)=\sum_{v\in \mathbb{L}_{\ee_n,A}}X(v)+\sum_{u\in A}X_0^*(u)-|\mathbb{O}_{\ee_n,A}|.
\end{equation}
\end{lemma}
\proof We will prove \eqref{e:AXs}  by reduction to $n$.
It is true for \eqref{e:AXs} when $n=0$ since $\TT_0^{\ee_n}=\{e_0\}$, $\Theta^{\{\ee_0\}}X(\ee_0)=X_0^*(\ee_0)$ and
 $\mathbb{O}_{\ee_0, \{\ee_0\}}=\emptyset=\mathbb{L}_{\ee_0,\{\ee_0\}}$.
 Let $\ell\ge 1$. Assume  that \eqref{e:AXs} holds true for  $n=\ell-1$, and we will prove it still holds true for $n=\ell$.

 Fix $i\ge 0$ and $A\subset \mathbb{T}_0^{\ee_{\ell}}$ with $A\not=\emptyset$. As before, set $[m]=\{1,2,\cdots,m\}$,  ~~$A_{j}= A\cap \mathbb{T}_0^{\ee_{\ell}^{(j)}}$ for $j\in [m]$ and $I=\{j\in [m]:   A_{j}\not=\emptyset\}$.
 Suppose event $\{\nabla^A 1_{\{X(\ee_{\ell})=i\}}\not=0\}$ occurs.

 By (2) of Lemma \ref{l:AXs00},  we have     $\Theta^A(\ee_{\ell})\ge 1$, which implies
\begin{align*}
  \Theta^A(\ee_{\ell})=\sum_{j=1}^m \Theta^A X(\ee_{\ell}^{(j)})-1
=\sum_{j\in I} \Theta^{A_j} X(\ee_{\ell}^{(j)})+\sum_{j\in [m]\backslash I}   X(\ee_{\ell}^{(j)})-1.
\end{align*}
By (3) of  Lemma \ref{l:AXs00}, there exist some integers $x_1\ge 0, \cdots, x_m\ge 0$  such that
\begin{equation}\label{e:Aj_not0}
\nabla^{A_{j}} 1_{\{ X(\ee_{\ell}^{(j)})=x_j\}}\not=0 {\rm ~for~all~}j\in [m].
\end{equation}
Note that $|\ee_{\ell}^{(j)}|=\ell-1$ and $A_j\subset \TT_0^{\ee_\ell^{(j)}}$ with $A_j\not=\emptyset$ for   $j\in I$. So by \eqref{e:Aj_not0} and
 the assumption that \eqref{e:AXs} holds true for $n=\ell-1$, we can get
$$
\Theta^{A_{j}} X(\ee_{\ell}^{(j)}) =\sum_{v\in \mathbb{L}_{\ee_{\ell}^{(j)},A_{j}}}X(v)+\sum_{u\in A_{j}}X_0^*(u)-|\mathbb{O}_{\ee_{\ell}^{(j)},A_{j}}|, \quad j\in I.
$$

 Combining these equalities together, we obtain that conditioned on $\{\nabla^A 1_{\{X(\ee_{\ell})=i\}}\not=0\}$,
\begin{align*}
\Theta^A(\ee_{\ell})=&\sum_{j\in I} (\sum_{v\in \mathbb{L}_{\ee_{\ell}^{(j)},A_{j}}}X(v)+\sum_{u\in A_{j}}X_0^*(u)-|\mathbb{O}_{\ee_{\ell}^{(j)},A_{j}}|)+\sum_{j\in [m]\backslash I} X(\ee_{\ell}^{(j)})-1.
\end{align*}
Since $A=\bigcup\limits_{j\in I} A_j$,  $\mathbb{O}_{\ee_{\ell},A}=\{e_{\ell}\}\cup \bigcup\limits_{j\in I}  \mathbb{O}_{\ee_{\ell}^{(j)}, A_{j}}$ and $\mathbb{L}_{\ee_{\ell},A}=\{e_{\ell}^{(j)}:  j\in [m]\setminus I \}   \cup \bigcup\limits_{j\in I}  \mathbb{L}_{\ee_{\ell}^{(j)}, A_{j}}$, we obtain further
\begin{align*}
\Theta^A(\ee_{\ell})= \sum_{v\in \mathbb{L}_{\ee_{\ell},A}}X(v)+\sum_{u\in A}X_0^*(u)-| \mathbb{O}_{\ee_{\ell},A}|.
\end{align*}
Such  we prove that \eqref{e:AXs} holds true for $n=\ell$,
and finish the proof of the lemma.
\qed\\

\quad\quad By  Lemmas \ref{l:AXs00} and \ref{l:AXs}, $\sum_{v\in \mathbb{L}_{\ee_n,A}}X(v)+\sum_{u\in A}X_0^*(u)-|\mathbb{O}_{\ee_n,A}|\ge 1$ conditioned on  $\{\nabla^A 1_{\{X(\ee_n)=0\}}\not=0\}$. So, we go ahead to study some quantity   related to those $|\mathbb{O}_{\ee_n,A}|$.
\begin{lemma}\label{l:branch_number} Let $k\ge 0$.  Then
\begin{equation}\label{e:branch_number}
\sum_{A\subset \mathbb{T}_0^{\ee_n}:|A|=k} m^{-|\mathbb{O}_{\ee_n,A}|}\le m^{k^m} n^{(k-1)^+}, \quad  n\ge 1.
\end{equation}
\end{lemma}
\proof We will prove \eqref{e:branch_number} by induction to $k$.  When $n\ge 1$ and  $A=\emptyset$, we have $\mathbb{O}_{\ee_n,A}=\emptyset$ and so \eqref{e:branch_number} holds true for    $k=0$.
When $n\ge 1$ and $v\subset \mathbb{T}_0^{\ee_n}$, we have $\mathbb{O}_{\ee_n,\{v\}}=\{v_1,\cdots,v_n\}$. So
$$
\sum_{A\subset \mathbb{T}_0^{\ee_n}:|A|=1} m^{-|\mathbb{O}_{\ee_n,A}|}=\sum_{v\in \mathbb{T}_0^{\ee_n}} m^{-|\mathbb{O}_{\ee_n,\{v\}}|}=\sum_{v\subset \mathbb{T}_0^{\ee_n}} m^{-n}=1.
$$
Hence \eqref{e:branch_number} holds true for   $k=1$.

Now let $\ell\ge 2$ and  assume that
\eqref{e:branch_number} holds true for each $k<\ell$, and we will show it still holds true for $k=\ell$.
Let  $n\ge 1$, $A\subset \mathbb{T}_0^{\ee_n}$ with $|A|=\ell$. Denote by $w$  the first common offspring of $A$, and set $A_j=A\cap \mathbb{T}_0^{w^{(j)}}$ for each $1\le j\le m$. Then
$$
\sum_{j=1}^m |A_j|=\ell,\quad |\mathbb{O}_{\ee_n, A}|=n-|w|+1+\sum_{j=1}^m|\mathbb{O}_{w^{(j)},A_j}|.
$$
Since $|A|=\ell\ge 2$, we  have  $1\le |w|\le n$, $\sum\limits_{j=1}^m (|A_j|-1)^+\le \ell-2$ and
 $\max\limits_{1\le j\le m}|A_j|\le\ell-1$.
So, writing $\Lambda_\ell=\{(\ell_1,\cdots,\ell_m):\sum\limits_{j=1}^m \ell_j=\ell, \sum\limits_{j=1}^m (\ell_j-1)^+\le \ell-2 {\rm~ and~} 0\le \ell_j\le\ell-1 {\rm ~for ~} j\}$, we get
\begin{align*}
\sum_{A\subset \mathbb{T}_0^{\ee_n}:|A|=\ell} m^{-\mathbb{O}_{\ee_n,A}}=&\sum_{h=1}^n~\sum_{w\in \mathbb{T}^{\ee_n}_h}~\sum_{(\ell_1,\cdots,\ell_m)\in \Lambda_\ell}~\sum_{ A_j\subset \mathbb{T}_0^{w^{(j)}} {\rm ~with~} |A_j|=\ell_j {\rm ~for~} j\le m} m^{-n+h-1- \sum_{j=1}^m|\mathbb{O}_{w^{(j)},A_j}|    }\\
=&\sum_{h=1}^n~\sum_{w\in \mathbb{T}^{\ee_n}_h}m^{-n+h-1}~\sum_{(\ell_1,\cdots,\ell_m)\in \Lambda_\ell}~\prod_{j=1}^m \left(\sum_{ A_j\subset \mathbb{T}_0^{w^{(j)}} {\rm ~with~} |A_j|=\ell_j} m^{-  |\mathbb{O}_{w^{(j)},A_j}|    }\right).
\end{align*}
Let $1\le j\le m$.  If $w\in \TT_1^{\ee_n}$, then $\mathbb{T}_0^{w^{(j)}}=\{w^{(j)}\}$ and $ \mathbb{O}_{w^{(j)},A_j}=\emptyset$, which implies
$$\sum_{ A_j\subset \mathbb{T}_0^{w^{(j)}} {\rm ~with~} |A_j|=\ell_j} m^{-  |\mathbb{O}_{w^{(j)},A_j}|    }=1_{\{\ell_j\le 1\}}\le m^{\ell_j^m} |w|^{(\ell_j-1)^+}.$$
While if $ w\in \TT^{\ee_n}\backslash(\TT_0\cup \TT_1)$ then
$  |w^{(j)}|\ge 1$. By $\ell_j\le \ell-1$,    we can use the assumption that  \eqref{e:branch_number} holds ture for  $k<\ell$ to get
$$
\sum_{ A_j\subset \mathbb{T}_0^{w^{(j)}} {\rm ~with~} |A_j|=\ell_j} m^{-  |\mathbb{O}_{w^{(j)},A_j}|    }\le m^{\ell_j^m} (|w^{(j)}|)^{(\ell_j-1)^+}\le m^{\ell_j^m} |w|^{(\ell_j-1)^+}.
$$
Therefore,
\begin{align*}
\sum_{A\subset \mathbb{T}_0^{\ee_n}:|A|=\ell} m^{-\mathbb{O}_{\ee_n,A}}\le &\sum_{h=1}^n~\sum_{w\in \mathbb{T}^{\ee_n}_h}m^{-n+h-1}~\sum_{(\ell_1,\cdots,\ell_m)\in \Lambda_\ell}~\prod_{j=1}^m  (m^{\ell_j^m} |w|^{(\ell_j-1)^+})\\
 \le &\sum_{h=1}^n~\sum_{w\in \mathbb{T}^{\ee_n}_h}m^{-n+h-1}~\sum_{(\ell_1,\cdots,\ell_m)\in \Lambda_\ell}~ m^{\sum_{j=1}^m\ell_j^m} n^{\sum_{j=1}^m(\ell_j-1)^+}.
\end{align*}
Since $(\ell_1,\cdots,\ell_m)\in \Lambda_\ell$, we have $\sum_{j=1}^m(\ell_j-1)^+\le \ell-2$ and $\sum_{j=1}^m \ell_j=\ell$.
Noting that $(a+b-1)^m+1=a^m+\sum_{r=1}^{m} \binom{m}{r} a^{m-r} (b-1)^{r}+1\ge  a^m+\sum_{r=1}^{m} \binom{m}{r}  (b-1)^{r}+1=a^m+b^m$ for any $a\ge 1$, $b\ge 1$ and $m\ge 2$, we get  $\sum_{j=1}^m \ell_j^m\le (\ell-1)^m+1\le \ell^m-\ell$. Hence
 \begin{align*}
   \sum_{A\subset \mathbb{T}_0^{\ee_n}:|A|=\ell} m^{-\mathbb{O}_{\ee_n,A}}\le&  \sum_{h=1}^n~|\mathbb{T}^{\ee_n}_h|~m^{-n+h-1}~|\Lambda_\ell|~ m^{\ell^m-\ell} n^{\ell-2}.
 \end{align*}
Hence by $| \mathbb{T}^{\ee_n}_h|=m^{n-h}$ and $|\Lambda_\ell|\le m^\ell$, we obtain
$$
\sum_{A\subset \mathbb{T}_0^{\ee_n}:|A|=\ell} m^{-\mathbb{O}_{\ee_n,A}}\le  m^{\ell^m-1} n^{\ell-1}.
$$
Such   \eqref{e:branch_number} holds true for $k=\ell$ and  we finish the proof.\qed\\

Now we are readily to prove Proposition \ref{t:pX1}. \\

{\it Proof of Proposition \ref{t:pX1}. }  Let $n\ge 1$ and  $k\ge 1$. To be conciseness,  for   $A\subset \mathbb{T}_0^{\ee_n}$ with $|A|=k$ we set
\begin{align*}
E_A:=&\{ \sum_{v\in \mathbb{L}_{\ee_n,A}}X(v)+\sum_{u\in A}X_0^*(u)-|\mathbb{O}_{\ee_n,A}|\ge 1\},\\
P_A:=& \P(X(\ee_n)=0, ~X(v)\le n-|v| {\rm ~for~any~} v\in \mathbb{L}_{\ee_n,A}, ~X|_A=0, ~E_A).
\end{align*}
Then by Lemmas \ref{l:Dev_F_n}, \, \ref{l:AXs00}  and
\ref{l:AXs},
\begin{align}
\nonumber|\frac{d^k}{dp^k}\P(X_n=0)|\le&  \frac{2^kk!}{(1-p)^k}\sum_{A\subset \mathbb{T}_0^{\ee_n}:~|A|=k}  \P(\nabla^{A} 1_{\{X_n=0\}}\not=0,~ X|_A=0)\\
\nonumber \le& \frac{2^kk!}{(1-p)^k}\sum_{A\subset \mathbb{T}_0^{\ee_n}:~|A|=k}  \P(X(v)\le n-|v| {\rm ~for~any~} v\in \mathbb{T}^{\ee_n},~ X|_A=0,~E_A)\\
 \label{d^k_sep} \le& \frac{2^k k!}{(1-p)^k}\sum_{A\subset \mathbb{T}_0^{\ee_n}:~|A|=k}  P_A.
\end{align}
Next, we will prove the statements  (1), (2) and (3) of Proposition \ref{t:pX1} separately.

(1) Since  $P_A\le \P(X(\ee_n)=0)=\P(X_n=0)$ for  $A\subset \mathbb{T}_0^{\ee_n}$ with $|A|=k$,
 $$|\frac{d^k}{dp^k}\P(X_n=0)|\le \frac{2^k k!}{(1-p)^k} \sum_{A\subset \mathbb{T}_0^{\ee_n}:~|A|=k} \P(X_n=0)\le \frac{2^k k!}{(1-p)^k}m^{kn}\P(X_n=0).$$

(2) Since   $X(\ee_n)\ge \sum_{v\in\mathbb{T}^{\ee_n}_0}X(v)-m^n$,   we get
\begin{align*}
P_A\le& \P(\sum_{v\in\mathbb{T}^{\ee_n}_0}X(v)-m^n\le 0,  X|_A=0)\\
=& \P(\sum_{v\in\mathbb{T}^{\ee_n}_0\backslash A}X(v)\le m^n)\prod_{v\in A}\P(X(v)=0)\\
=& \P(\sum_{i=1}^{m^n-k} X_{0,i}\le m^n) ~(1-p)^k,
\end{align*}
 where $X_{0,i}, i\ge 1$ are i.i.d copies of $X_0$.
 Taking  the above inequality with \eqref{d^k_sep}  and $|\{A\subset \mathbb{T}_0^{\ee_n}:~|A|=k\}|\le m^{nk}$ ,  we obtain (2).

(3) Fix $\delta\in (0,\frac{1}{2})$ and $c_2:=\E(m^{(1-\delta)X_0^*})<\infty$.   Choose
$$
n_1:=\min\{n\ge 10: 2 k c_2 \,  m^{\delta^{-4}(1+\ln n)^2 (m-1)} \, m^{k^{m-1}} n\le   m^{\delta  n}\}.
$$
Let $n\ge n_1$. Since    $m^{(1-2\delta)(\sum_{v\in \mathbb{L}_{\ee_n,A}}X(v)+\sum_{u\in A}X_0^*(u)-|\mathbb{O}_{n,A}|)}\ge 1_{E_A}$,
$$
 P_A\le \E\left(m^{(1-2\delta)(\sum_{v\in \mathbb{L}_{\ee_n,A}}X(v)+\sum_{u\in A}X_0^*(u)-|\mathbb{O}_{\ee_n,A}|)}\prod_{u\in A}1_{\{X(u)=0\}}\prod_{v\in\mathbb{L}_{\ee_n,A}}1_{\{   X(v)\le n-|v| \}}\right).
$$
 Since  these   $X_0^*(u), 1_{\{X(u)=0\}}, X(v) $ for $u\in A$ and $v\in \mathbb{L}_{\ee_n,A}$ are independent,
\begin{align*}
  P_A\le&     \prod_{u\in A}[\E(m^{(1-2\delta)X_0^*(u)})\P(X(u)=0)] ~\cdot~
  \prod_{v\in  \mathbb{L}_{\ee_n,A}} [\E(m^{(1-2\delta)X(v)}1_{\{ X(v)\le  n-|v|\}})] ~\cdot~ m^{(2\delta-1)|\mathbb{O}_{\ee_n,A}|}\\
  =&    [\E(m^{(1-2\delta)X_0^*})(1-p)]^{|A|}~~
   \prod_{i=0}^{n-1} [\E(m^{(1-2\delta)X_i}1_{\{ X_i\le  n-i\}})]^{|\mathbb{L}_{\ee_n,A}\cap \TT_i|} ~~ m^{(2\delta-1)|\mathbb{O}_{\ee_n,A}|}.
\end{align*}
Since  $\E(m^{(1-2\delta)X_0^*})\le \E(m^{(1-\delta)X_0^*})= c_2$, ~$|A|=k$,  $|\mathbb{O}_{\ee_n,A}|\le |A|n$   and
  $|  \mathbb{L}_{\ee_n,A}\cap \TT_i|\le (m-1)|A|$ for $i<n$, the above inequality can be simplified as
\begin{align}\label{e:Prod1-2delta}
 P_A\le& [c_2 (1-p)]^k \left(\prod_{i=0}^{n-1} \E(m^{(1-2\delta)X_i}1_{\{ X_i\le n-i\}})\right)^{(m-1)k} ~m^{2\delta nk-|\mathbb{O}_{\ee_n,A}|}.
\end{align}

So we need estimate  $\prod_{i=0}^{n-1} \E(m^{(1-2\delta)X_i}1_{\{ X_i\le n-i\}})$. Write  $M :=\lfloor n-\delta^{-1}\ln n\rfloor$  and $X_i^{(M)} :=X_i\wedge (M-i)$ for each $i<M$.  Then   for
 $i\le \lfloor n-\delta^{-2}(1+\ln n)\rfloor$,
\begin{align*}
 (1-\delta)X_i^{(M)}-(1-2\delta) X_i1_{\{ X_i\le n-i\}}\ge& (1-\delta)(M-i)-(1-2\delta)(n-i)\\
 =&-(1-\delta)(n-M)+\delta(n-i)\\
 \ge& -(1+\delta^{-1}\ln n)+\delta\cdot \delta^{-2}(1+\ln n)\\
 =&\delta^{-1}-1>0.
\end{align*}
It follows $\E(m^{(1-2\delta)X_i}1_{\{ X_i\le n-i\}})\le \E(m^{(1-\delta)X_i^{(M)}})$ for  $i\le \lfloor n-\delta^{-2}(1+\ln n)\rfloor$.
On the  other hand,   for   $\lfloor n-\delta^{-2}(1+\ln n)\rfloor<i<n$
we have  $m^{(1-2\delta)X_i}1_{\{X_i\le n-i\}}\le m^{\delta^{-2}(1+\ln n)}$. Hence
\begin{align*}
\prod_{i=0}^{n-1} \E(m^{(1-2\delta)X_i}1_{\{ X_i\le n-i\}})\le& \prod\limits_{i=\lfloor n-\delta^{-2}\ln n\rfloor+1}^{n-1}  m^{\delta^{-2}(1+\ln n)}~\prod\limits_{i=0}^{\lfloor n-\delta^{-2}(1+\ln n)\rfloor}  \E(m^{(1-\delta)X_i^{(M)}}) \cdot \\
\le& m^{\delta^{-4}(1+\ln n)^2}\prod\limits_{i=0}^{M-1}  \E(m^{(1-\delta)X_i^{(M)}}).
\end{align*}
Hence by \eqref{e:Prod1-2delta}, we have
\begin{align*}
P_A
\le&   [ c_2(1-p) ]^k \left(m^{\delta^{-4}(1+\ln n)^2} \prod_{i=0}^{M-1} \E(m^{(1-\delta)X_i^{(M)}})\right)^{(m-1)k}~m^{2\delta nk-|\mathbb{O}_{\ee_n,A}|}.
\end{align*}

By Lemma  \ref{l:branch_number}, we have  $\sum_{A\subset \mathbb{T}_0^{\ee_n}:|A|=k} m^{-\mathbb{O}_{\ee_n,A}}\le m^{k^m} n^{(k-1)^+}$. So, by  \eqref{d^k_sep},   \begin{align*}
|\frac{d^k}{dp^k}\P(X_n=0)|\le &  \frac{2^k k!}{(1-p)^k}\sum_{A\subset \mathbb{T}_0^{\ee_n}:~|A|=k}  P_A\\
\le& 2^k k! c_2^k \left(  m^{\delta^{-4}(1+\ln n)^2 }\prod_{i=0}^{M-1} \E(m^{(1-\delta)X_i^{(M)}})\right)^{(m-1)k} \sum_{A\subset \mathbb{T}_0^{\ee_n}:|A|=k} m^{2\delta nk-|\mathbb{O}_{\ee_n,A}|} \\
\le&   [2 k c_2    m^{\delta^{-4}(1+\ln n)^2 (m-1)}]^k \left(\prod_{i=0}^{M-1} \E(m^{(1-\delta)X_i^{(M)}})\right)^{(m-1)k} ~m^{2\delta nk}m^{k^m} n^{(k-1)^+}.
\end{align*}
By the definition of $n_1$, we then prove that (3) holds true for $n\ge n_1$. We have completed the proof of the proposition.\qed\\



\section{Proof of Proposition \ref{t:pX2}}
  Proposition \ref{t:pX2}, crudely speaking, tells us   $\prod\limits_{i=0}^{M-1}\E(~m^{(1-\delta)(X_i\wedge (M-i))}~)^{m-1}$ growing to infinity with speed at most $m^{o(M)}$   uniformly in $p\in (0,1)$ as $M\rightarrow\infty$ or   $\P(X_n=0)$  decaying  with  at least doubly exponential speed eventually.

The  idea of estimating  the upper bounds of $\prod\limits_{i=0}^{M-1}\E(~m^{(1-\delta)(X_i\wedge (M-i))}~)^{m-1}$ comes from \cite{collet-eckmann-glaser-martin}, \cite{bmxyz_questions} and  \cite{CDDHLS2021}. Inequality  $\prod\limits_{i=0}^{M-1}\E(m^{X_i})^{m-1}\le cM^2$ was proved to be true for all  $M\ge 1$ if    system $\{X_n, n\ge 0\}$ is  subcritical or  critical, see
   Collet, Eckmann, Glaser and Martin \cite{collet-eckmann-glaser-martin} and  Chen,   Derrida, Hu, Lifshit and Shi \cite{bmxyz_questions}.
 Recently in \cite{CDDHLS2021}, Chen, Dagard, Derrida, Hu, Lifshit and Shi showed that when $p\rightarrow p_c+$, it keeps $\prod\limits_{i=0}^{M-1}\E(m^{X_i})^{m-1}\le c'M^2$ for   $1\le M\le c''(p-p_c)^{-1/2}$.
 However, we need a uniform upper bounds in $p\in(0,1)$ as $M\rightarrow\infty$ in the current paper. We have to exclude some situation; Luckily,  this situation play a role  under which   $\P(X_n=0)$  decays  with  at least doubly exponential speed eventually, see  Theorem  \ref{l:trucation}.

To obtain the decay speed of   $\P(X_n=0)$ as required,
we collect some properties for  general Derrida-Retaux systems;
 See Fact \ref{fact1}, Lemma  \ref{l:zero1} and Lemma \ref{l:M+zero} below. There, we  use $K_i, 1\le i\le 6$ to stand for some  constants which are independent of $n$.

\begin{fact}\label{fact1}(\cite[Theorem 6.5 and Lemma 2.4]{CDDHLS2021})  Fix $0<\alpha< \beta<1, \gamma>0$ and $\eta\in (0,\frac{1}{3(m-1)}]$.
Let $\{X_n, n\ge 0\}$ and $\{Y_n, n\ge 0\}$ be two Derrida-Retaux systems with $(m-1)\E(Y_0m^{Y_0})=\E(m^{Y_0})$ and  $\P(X_0=k)\ge \P(Y_0=k)$ for all $k\ge 1$.    Assume that     $\P(Y_0=0)\in [\alpha,\beta]$,  $\E(Y_0^3m^{Y_0})\le \gamma$ and $\E(X_0-Y_0)\ge \eta$.  Then
there exists some constant  $K_1=K_1(m,\alpha,\beta,\gamma,\eta)\in \mathbf{Z}^+$ such that
$$
\max\limits_{0\le j\le K_1} \E(X_j)\ge 2.
$$
\end{fact}

\begin{lemma}\label{l:zero1} Fix $\lambda\in \{2,3,\cdots\}$ and  $\theta>0$. Let $\{X_n, n\ge 0\}$   be a Derrida-Retaux system  with  $X_0\le \lambda$ and $\E(((m-1)X_0-1)m^{X_0})\ge \theta$.
Then there exist  some constants  $K_2=K_2(m,\lambda,\theta)\in \mathbf{Z}^+$ and    $K_3=K_3(m,\lambda,\theta)>0$  such that
$$
\P(X_{K_2+n}=0)\le   e^{-K_3 m^n}, \quad n\ge 0.
$$
\end{lemma}
\proof Assume   $X_0\le \lambda$ and $(m-1)\E(X_0m^{X_0})-\E(m^{X_0})\ge \theta>0$. Then
$$
\P(X_0\ge 1)\ge \frac{(m-1)\E(X_0 m^{X_0})}{(m-1)\lambda m^{\lambda}}\ge \frac{\E(m^{X_0})+\theta}{(m-1)\lambda m^{\lambda}}\ge \frac{1}{(m-1)\lambda m^{\lambda}}.
$$
 Let $Z_0$ be a Bernoulli random variable which  satisfies $\P(Z_0=0)=  \frac{(m-1)\E(X_0m^{X_0})-\E(m^{X_0})}{1+(m-1)\E(X_0m^{X_0})-\E(m^{X_0})  }$ and   $\P(Z_0=1)=\frac{1}{1+(m-1)\E(X_0m^{X_0})-\E(m^{X_0})  }$.
Then
$$
\P(Z_0=0)\ge \frac{\theta}{1+\theta} \quad {\rm and}\quad   \P(Z_0=1)\ge  \frac{1}{1+(m-1)\lambda m^\lambda}.
$$
Assume further $Z_0$ is independent of $X_0$ and let $\{Y_n, n\ge 0\}$ be the Derrida-Reatux system with   $Y_0:=X_0Z_0$.  Then system $\{Y_n, n\ge 0\}$ is  critical, this is because
\begin{align*}
 \E(((m-1)Y_0-1)m^{Y_0})  =\E(((m-1)X_0-1)m^{X_0})\P(Z_0=1)- \P(Z_0=0)  =0.
 \end{align*}

 We want to
 apply Fact \ref{fact1} for  systems $\{X_n, n\ge 0\}$ and $\{Y_n, n\ge 0\}$, so that  we should check those conditions.
  By construction, we have the following:
\begin{align*}
\P(Y_0=k)=&\P(Z_0=1)\P(X_0=k)\le \P(X_0=k) {\rm ~~for~~}k\ge 1;\\
\P(Y_0=0)\ge &\P(Z_0=0)\ge \frac{\theta}{1+\theta};\\
 \P(Y_0=0)=&1-\P(X_0\ge 1)\P(Z_0=1)\le 1- \frac{1}{(m-1)\lambda m^{\lambda} \, (1+(m-1)\lambda m^\lambda)  };\\
 \E(Y_0^3 m^{Y_0})\le& \E(X_0^3 m^{X_0})\le \lambda^3 m^\lambda; {\rm ~and~}\\
 \E(X_0-Y_0)\ge&\E(1_{\{X_0\ge 1, Z_0=0\}})= \P(X_0\ge 1)\P(Z_0=0)\ge \frac{\theta}{(m-1)\lambda m^{\lambda} \, (1+\theta)}.
\end{align*}
So, we can apply  Fact \ref{fact1} with $\alpha=\frac{\theta}{1+\theta}, \, \beta= 1- \frac{1}{(m-1)\lambda m^{\lambda}(1+(m-1)\lambda m^\lambda)  }, \, \gamma=\lambda^3m^\lambda$ and $\eta=  \frac{\theta}{(m-1)\lambda m^{\lambda}(1+\theta)}$  to get
$$
\max_{0\le j\le K_1(m,\alpha,\beta,\gamma,\eta)}\E(X_{n})\ge 2.
$$
We choose  $K_2 :=K_1(m,\alpha,\beta,\gamma,\eta)$, so that the value of  $K_2$ depends only on $(m,\lambda,\theta)$.  Since  $m\ge 2$ and $\E(X_{n+1})\ge m\E(X_n)-1$ for all $n$, we have $\E(X_{K_2})\ge 2$ always.

 Since $0\le X_0\le \lambda$, we have $0\le X_{K_2}\le m^{K_2} \lambda$. Let $n\ge 0$.
Since $X_{K_2+n}$ is stochastically greater than $\sum_{i=1}^{m^n} X_{K_2,i}-m^n$,
$$
\P(X_{K_2+n}=0)\le \P(\sum_{i=1}^{m^n} X_{K_2,i}-m^n\le 0),
$$
where  $X_{K_2,i}, \, i\ge 1$ are i.i.d copies of $X_{K_2}$. By   $\E(X_{K_2})\ge 2$ and $0\le X_{K_2}\le m^{K_2}\lambda$, we apply  the   Hoeffding's inequality to get
$$
\P(\sum_{i=1}^{m^n} X_{K_2,i}\le  m^n)\le \P(\sum_{i=1}^{m^n} X_{K_2,i}\le \sum_{i=1}^{m^n} \E(X_{K_2,i})- m^n) \le   e^{-\frac{2m^n}{(m^{K_2}\lambda)^2}}.
$$
Taking $K_3 :=\frac{2}{m^{2K_2}\lambda^2}$, we draw out the conclusion of the lemma.\qed\\

\begin{lemma}\label{l:M+zero}  Fix $\delta\in (0,1)$ and  $\theta>0$. Let $\{X_n, n\ge 0\}$   be a Derrida-Retaux system. Assume that
  $\E(((m-1)X_0-1)s^{X_0})\ge 0$,    $\E(((m-1)X_0-1)s^{X_0}1_{\{X_0\ge 1\}})\ge \theta$ and $X_0\le M$ for some
$s\in [1, m-\delta]$  and some $M\in \{2,3,\cdots\}$. Then there exist some constants $K_4=K_4(m,\delta,\theta)\in \mathbf{Z}^+$ and   $K_5=K_5(m,\delta,\theta)>0$  such that
$$
\P(X_{M+K_4+n}=0)\le   e^{- K_5 m^{ n}}, \quad n\ge 0.
$$
\end{lemma}
\proof Fix  $s\in [1, m-\delta]$ and    $M\in \{2,3,\cdots\}$ which satisfy   $\E(((m-1)X_0-1)s^{X_0})\ge 0$,   $\E(((m-1)X_0-1)s^{X_0}1_{\{X_0\ge 1\}})\ge \theta$ and $X_0\le M$.
  Since $X_0\le M$, we need only  consider  Case I  ~(\, $\max\limits_{1\le k\le M} m^k\P(X_0=k+2)\ge  1$ \, ) and Case II ~(\, $\sup\limits_{k\ge 1} m^k\P(X_0=k+2)\le 1$\, ).

Case I: Suppose there exists  some $1\le k\le M$ such that   $ \P(X_0=k+2)\ge  m^{-k}$. Since $X_k$ is stochastically greater than $\max_{1\le i\le m^k} X_{0,i}-k$, where $X_{0, i}, ~i\ge 1$ are i.i.d copies of $X_0$,
  \begin{align*}
 \P(X_k\ge 2)\ge&  \P(X_{0,i}=k+2 {\rm ~for~ some~}1\le i\le m^k)\\
 =&  1-(1-\P(X_0=k+2))^{m^k} \\
 \ge &  1-(1-m^{-k})^{m^k}.
 \end{align*}
 Since $\sup\limits_{x\in (0,1)}(1-x)^{x^{-1}}=\lim\limits_{x\rightarrow 0+}(1-x)^{x^{-1}}=e^{-1}$,      we  have $\P(X_k\ge 2)\ge 1-e^{-1}$.
 Let $\{\tilde X_n, n\ge 0\}$  be the  Derrida-Retaux system  with $\tilde X_0:=X_k\wedge 2$.  Then  $\tilde X_0\le 2$ and
\begin{align*}
\E(((m-1)\tilde X_0-1)m^{\tilde X_0})\ge& (2(m-1)-1)m^2\P(X_k\ge 2)-\P(X_k=0)\\
\ge& 2^2(1-e^{-1})-1\ge  1.
\end{align*}
 So,  we can apply
  Lemma \ref{l:zero1} for  system $\{\tilde X_n, n\ge 0\}$  to get
  $$
  \P(\tilde X_{K_2(m,2,1)+n}=0)\le e^{-K_3(m,2,1) m^n}, \quad n\ge 0.
  $$
  Since $\tilde X_0\le X_k$, we have $\tilde X_n\le X_{k+n}$ for all $n\ge 0$. So,
  $$
  \P( X_{k+K_2(m,2,1)+n}=0)\le \P(\tilde X_{K_2(m,2,1)+n}=0) \le e^{-K_3(m,2,1) m^n}.
  $$
  Since $k\le M$, the above inequality can be rewrote as
  $$
   \P( X_{M+K_2(m,2,1)+n}=0)\le e^{-K_3(m,2,1) m^{n+M-k}}\le e^{-K_3(m,2,1) m^{n}}.
  $$

Case II: Suppose    $\P(X_0=k+2)\le  m^{-k}$ for $k\ge 1$.   Write $t=m-\frac{1}{2}\delta$  and set
 $$K_6:=\min\{\ell\ge 2: \sum\limits_{k=\ell+1}^\infty ((m-1)k-1)(\frac{t}{m})^k \le \frac{\theta\delta}{4m^3}\}.$$
 Then
 \begin{align*}
 \E(((m-1)X_0-1)t^{X_0}1_{\{X_0\ge K_6+1\}})=&\sum\limits_{k=K_6+1}^\infty ((m-1)k-1)t^k\P(X_0=k)\\
 \le& \sum\limits_{k=K_6+1}^\infty ((m-1)k-1)t^km^{-k+2}\\
 \le&    \frac{\theta\delta}{4m}.
 \end{align*}
 On the other hand,  $\frac{t}{s}\ge \frac{m-\frac{\delta}{2}}{m-\delta}\ge 1+\frac{\delta}{2m}$  since $s\in [1, m-\delta]$. By
 $\E(((m-1)X_0-1)s^{X_0})\ge 0$ and
 $\E(((m-1)X_0-1)s^{X_0}1_{\{X_0\ge 1\}})\ge \theta$, we get
\begin{align*}
&\E(((m-1)X_0-1)t^{X_0})\\
=&\E(((m-1)X_0-1)(t^{X_0}-s^{X_0})1_{\{X_0\ge 1\}}) + \E(((m-1)X_0-1)s^{X_0})\\
\ge& \frac{\delta }{2m}\E(((m-1)X_0-1) s^{X_0}1_{\{X_0\ge 1\}})\\
\ge& \frac{\theta\delta}{2m}.
\end{align*}
Hence
\begin{align*}
\E(((m-1)X_0-1)t^{X_0}1_{\{0\le X_0\le K_6\}})\ge \frac{\theta\delta}{2m}-\frac{\theta\delta}{4m}=  \frac{\theta\delta}{4m}.
\end{align*}
Let  $\{ \hat{X}_n, n\ge 0\}$ be the Derrida-Retaux system  with $ \hat{X}_0:=X_0\wedge K_6$. Then $\hat{X}_0\le K_6$ and
$$
\E(((m-1)X_0-1)m^{\hat{X}_0})\ge \E(((m-1)X_0-1)t^{X_0}1_{\{0\le X_0\le K_6\}})\ge \frac{\theta\delta}{4m}.
$$
So we can apply Lemma \ref{l:zero1}  for  system  $\{ \hat{X}_n, n\ge 0\}$ to get
  $$
  \P( \hat X_{K_2(m,K_6,\frac{\theta\delta}{4m})+n}=0)\le e^{-K_3(m,K_6,\frac{\theta\delta}{4m}) m^n}, \quad n\ge 0.
  $$
  Since $\hat{X}_0\le X_0$ we have $\hat{X}_n\le X_n$ for all $n\ge 0$ which implies
  $$
    \P( X_{M+K_2(m,K_6,\frac{\theta\delta}{4m})+n}=0)\le \P( \hat X_{M+K_2(m,K_6,\frac{\theta\delta}{4m})+n}=0)\le e^{-K_3(m,K_6,\frac{\theta\delta}{4m}) m^n}.
  $$

  We have obtained the estimates for both cases. Taking $K_4 :=K_2(m,2,1)\vee K_2(m,K_6,\frac{\theta\delta}{4m})$ and $K_5 :=K_3(m,2,1)\wedge K_3(m,K_6,\frac{\theta\delta}{4m})$, we finish the proof of the lemma.
\qed\\


Now we return back to our setting.
 Recall that  $c_1:=\P(X_0^*\ge 2)>0$. When $m=2$, the Derrida-Retaux system has a fix point $\delta_1$, whose support is  focused on set $\{1\}$. We have to   cope with the special case. Set $n_2:=\lfloor \frac{\log(\frac{1}{\P(X_0^*\ge 2)})}{\log(\frac{5}{4})}\rfloor+1$.
\begin{lemma}\label{l:PX_n=1}   If $m=2$, then $\P(X_n=1)\le \frac{1}{2}$  for any $  n\ge n_2$.
\end{lemma}
\proof Fix $m=2$, then
\begin{align}
\nonumber\P(X_{n+1}=1)=&2\P(X_n=2)\P(X_n=0)+\P(X_n=1)^2\\
 \nonumber\le&\frac{1}{2}(\P(X_n=2)+\P(X_n=0))^2+\P(X_n=1)^2\\
\label{e:X_n=1}\le& \frac{1}{2}(1-\P(X_n=1))^2+\P(X_n=1)^2.
\end{align}

  Set $\ell:=\inf\{n\ge 0: \P(X_n=1)\le \frac{1}{2} \}$.
Then $\P(X_n=1)>\frac{1}{2}$  for  $0\le n<\ell$. So, by \eqref{e:X_n=1}
 for  $0\le n<\ell$ we have
\begin{align*}
1-\P(X_{n+1}=1)\ge&1-(\frac{1}{2}(1-\P(X_n=1))^2+\P(X_n=1)^2)\\
=& \frac{1}{2} (1+3\P(X_n=1))(1-\P(X_n=1))\\
\ge& \frac{5}{4}(1-\P(X_n=1)).
\end{align*}
It follows immediately,
$$
1-\P(X_{n+1}=1)\ge (\frac{5}{4})^{n+1}(1-\P(X_0=1)), \quad 0\le n< \ell.
$$
Since $\P(X_0^*\ge 2)>0$, we have $ 1-\P(X_0=1)\ge 1-\P(X_0^*=1)\ge   \P(X_0^*\ge 2)>0$. Hence
  $\ell$ is finite and satisfies
$$
1-\P(X_\ell=1)\ge (\frac{5}{4})^\ell (1-\P(X_0=1))\ge(\frac{5}{4})^\ell \P(X_0^*\ge 2).
$$
So,
$$
\ell\le \frac{\log(\frac{1}{\P(X_0^*\ge 2)})}{\log(\frac{5}{4})}\le n_2.
$$

On the other hand, by the definition of $\ell$, we have $\P(X_\ell=1)\le \frac{1}{2}$. Using \eqref{e:X_n=1} and the inequality $\frac{1}{2}(1-x)^2+x^2\le \frac{1}{2}$ for  $x\in [0,\frac{1}{2}]$, we show that   $\P(X_{n+1}=1)\le \frac{1}{2}$ once $\P(X_n=1)\le \frac{1}{2}$. Therefore, $\P(X_n=1)\le \frac{1}{2}$ for   $n\ge \ell$.
 We have completed the proof of the lemma.\qed\\

  We do not directly apply
 Lemma \ref{l:M+zero}
for system $\{X_n, n\ge 0\}$ since that $X_n$ is unbounded in general.  We need a  truncation.
Set
$$
X_i^{(M)} :=X_i\wedge (M-i), \quad M>i\ge 0.
$$

\begin{theorem}\label{l:trucation} Fix $\delta\in (0,1)$ and let $M\ge n_2+2$. Assume     $\E([(m-1)X_i^{(M)}-1]s^{X_i^{(M)}})\ge 0$ for some $i\in [n_2, M-2]\cap{\mathbf{Z}}$ and some $s\in [1,m-\delta]$. Then there exist constants $c_i=c_i(m,\delta)>0$, $i\in \{6,7\}$  such that
\begin{equation}\label{e:trucation11}
\P(X_n=0)\le c_7 e^{-c_6 m^{n-M}},\quad n\ge M.
\end{equation}
\end{theorem}
\proof Fix $M\ge n_2+2$,  $n_2\le i\le  M-2$ and   $s\in [1,m-\delta]$ with $\E([(m-1)X_i^{(M)}-1]s^{X_i^{(M)}})\ge 0$. Then it is true for the following statement   whose proof will be given  a little  later:  \begin{equation}\label{e:-theta+}
\E([(m-1)X_i^{(M)}-1]s^{X_i^{(M)}}1_{\{X_i^{(M)}\ge 1\}})\ge \frac{1}{4}.
\end{equation}
By admitting \eqref{e:-theta+}, we can apply  Lemma   \ref{l:M+zero} for the Derrida-Retaux system $(\tilde X_n, n\ge 0)$ with $\tilde X_0:=X_i^{(M)}$ to  get
 $$
 \P(\tilde X_{(M-i)+K_4(m,\delta,\frac{1}{4})+n}=0)\le e^{-K_5(m,\delta,\frac{1}{4}) m^n}, \quad n\ge 0.
 $$
 Since  $\tilde X_0\le X_i$, we have $\tilde X_n\le X_{n+i}$ for all $n\ge 0$. It follows immediately for any $n\ge M+K_4(m,\delta,\frac{1}{4})$,
 $$
 \P( X_{n}=0)\le  \P( \tilde X_{n-i}=0)\le e^{-K_5(m,\delta,\frac{1}{4}) m^{n-M-K_4(m,\delta,\frac{1}{4})}}.
 $$
 Set $c_6 :=K_5(m,\delta,\frac{1}{4}) m^{-K_4(m,\delta,\frac{1}{4})}$ and $c_7 :=e^{K_5(m,\delta,\frac{1}{4})}$.  Using the above inequality  and the fact $\P(X_n=0)\le 1$ for all $n$,   we draw out the conclusion of \eqref{e:trucation11}.

We are left to prove   \eqref{e:-theta+}.   If   $m\ge 3$,  then $((m-1)X_i^{(M)}-1)1_{\{X_i^{(M)}\ge 1\}}\ge \frac{1}{2}(m-1)X_i^{(M)}$.  By  $\E(((m-1)X_i^{(M)}-1)s^{X_i^{(M)}})\ge 0$ and $s\ge 1$, we then have
  \begin{align*}
  \E(((m-1)X_i^{(M)}-1)s^{X_i^{(M)}}1_{\{X_i^{(M)}\ge 1\}})\ge& \frac{1}{2} \E((m-1)X_i^{(M)}s^{X_i^{(M)}})  \ge \frac{1}{2}\E(s^{X_i^{(M)}})\ge \frac{1}{2}.
  \end{align*}
Now let $m=2$. Since $i\ge n_2$, we apply Lemma \ref{l:PX_n=1} to get $\P(X_i=1)\le \frac{1}{2}$. Since $M-i\ge 2$,
$$
\P(X_i^{(M)}=1)=\P(X_i=1)\le \frac{1}{2}.
$$
Since $s\ge 1$, we have $\E(((m-1)X_i^{(M)}-1)s^{X_i^{(M)}}1_{\{X_i^{(M)}\ge 1\}})\ge \P(X_i^{(M)}\ge 2)$. Since $\E(((m-1)X_i^{(M)}-1)s^{X_i^{(M)}})\ge 0$, we have $\E(((m-1)X_i^{(M)}-1)s^{X_i^{(M)}}1_{\{X_i^{(M)}\ge 1\}})\ge\P(X_i^{(M)}=0)$. Consequently,
  \begin{align*}
  \E(((m-1)X_i^{(M)}-1)s^{X_i^{(M)}}1_{\{X_i^{(M)}\ge 1\}})\ge& \frac{1}{2}((X_i^{(M)}=0)+\P(X_i^{(M)}\ge 2))\\
  =&\frac{1}{2}(1-\P(X_i^{(M)}=1))\ge \frac{1}{4}.
  \end{align*}
 Such  \eqref{e:-theta+} holds true for any $m\ge 2$, and we complete the proof of the theorem.  \qed\\

We will
 use  the method developed by  Collet, Eckmann, Glaser and Martin \cite{collet-eckmann-glaser-martin} to
 obtain an upper bound of $\prod\limits_{i=0}^{M-1} \E(m^{(1-\delta)X_i^{(M)}})$. Write $H_i^{(M)}(s) :=\E(s^{X_i^{(M)}})$ for $0\le i<M$ and $s\ge 0$.  As in  \cite{collet-eckmann-glaser-martin},   set
 \begin{align}\label{def:Delta}
 \Delta_i^{(M)}(s) :=&[H_i^{(M)}(s)-s(s-1)H_i^{(M)'}(s)]\\
 \nonumber&~~-\frac{(m-1)(m-s)}{m}[2sH_i^{(M)'}(s)+s^2H_i^{(M)''}(s)].
 \end{align}

\begin{lemma}\label{l:Cauchy_X} Fix $\delta\in (0,\frac{1}{16m})$ and $s=m^{1-\delta}$.   Let $n_2\le i\le M-2$ and set
  $ s_i>0$ for the value which satisfies
 $ \E([(m-1)X_i^{(M)}-1] s_i^{X_i^{(M)}})=0$. If  $s_i\ge m-m\delta^3$,  then
\begin{equation}\label{e:Cauchy_Scharz_D}
[H_i^{(M)}(s)-(m-1)sH_i^{(M)'}(s)]^2\le 2H_i^{(M)}(0)\Delta_i^{(M)}(s);
\end{equation}
\begin{equation}\label{e:Delta_low_bound}
\Delta_i^{(M)}(s)\ge \frac{\delta^2}{128}.
\end{equation}
\end{lemma}
\proof  Fix $\delta\in (0,\frac{1}{16m}), s=m^{1-\delta}$ and $n_2\le i\le M-2$ with $s_i\ge m-m\delta^3$. Write $x_i=\frac{s}{s_i}$ for conciseness.
 Then
$$
x_i\le \frac{m^{1-\delta}}{m-m\delta^3}=\frac{m^{-\delta}}{1-\delta^3}\le \frac{2^{-\delta}}{1-\delta^3}\le 1-\frac{\delta}{4}.
$$
  Since $\E(((m-1)X_i^{(M)}-1)s_i^{X_i^{(M)}})=0$,
\begin{align*}
 H_i^{(M)}(s)=&(m-1)s_iH_i^{(M)'}(s_i)- H_i^{(M)}(s_i)+H_i^{(M)}(s)\\
=&\sum_{k\ge 1}(km-k-1+x_i^k)s_i^k\P(X_i^{(M)}=k).
\end{align*}
Hence
$$H_i^{(M)}(0)=\sum\limits_{k\ge 1}(km-k-1)s_i^k\P(X_i^{(M)}=k),
$$
$$
H_i^{(M)}(s)-(m-1)sH_i^{(M)'}(s)=\sum_{k\ge 1}(1-x_i^k) (km-k-1)s_i^k\P(X_i^{(M)}=k)
$$
and
$$
\Delta_i^{(M)}(s)=\sum_{k\ge 1} (1-(k+1)x_i^k+\frac{s_i}{m}kx_i^{k+1})(km-k-1)s_i^k\P(X_i^{(M)}=k).
$$
By $s_i\ge m-m\delta^3$, we get
\begin{align}
\label{e:CauS3}\Delta_i^{(M)}(s)\ge &\sum_{k\ge 1} (1-(k+1)x_i^k+kx_i^{k+1}-\delta^3 k x_i^{k+1}) (km-k-1)s_i^k\P(X_i^{(M)}=k).
\end{align}
Set $\eta:=\sup_{k\ge 1} \frac{(1-x_i^k)^2}{ 1-(k+1)x_i^k+kx_i^{k+1}- \delta^3 k x_i^{k+1} }$. Then by the Cauchy inequality,
$$
[H_i^{(M)}(s)-(m-1)sH_i^{(M)'}(s)]^2\le \eta H_i^{(M)}(0)\Delta_i^{(M)}(s).
$$

We need  an upper bound of    $\eta$.  An elementary calculation gives for   $k\ge 1$ and $x\in {\mathbf{R}}$,
  $$2(1-(k+1)x^k+kx^{k+1})-(1-x^k)^2=k(1-x)^2x^{k-1}+(1-x) \sum_{\ell=0}^{k-1} (x^\ell-x^k)(1-x^{k-1-\ell}).$$
Since $0\le x_i\le 1- \frac{1}{4}\delta$ and $\delta<\frac{1}{16m}\le \frac{1}{32}$, we have
 \begin{align*}
 2(1-(k+1)x_i^k+kx_i^{k+1})-(1-x_i^k)^2\ge k(1-x_i)^2x_i^{k-1}\ge k (\frac{\delta}{4})^2 x_i^{k+1}\ge 2\delta^3k x_i^{k+1}.
 \end{align*}
So,
\begin{equation}\label{e:eta=2}
(1-x_i^k)^2\le 2(1-(k+1)x_i^k+ kx_i^{k+1} -\delta^3 kx_i^{k+1} ), \quad k\ge 1,
\end{equation}
which implies $\eta\le 2$. So we have  \eqref{e:Cauchy_Scharz_D}.

We are left to prove \eqref{e:Delta_low_bound}. By \eqref{e:CauS3} and \eqref{e:eta=2} we have
\begin{align*}
\Delta_i^{(M)}(s)\ge& \frac{1}{2} \sum_{k\ge 1} (1-x_i^k)^2   (km-k-1) s_i^k \P(X_i^{(M)}=k)\\
\ge&\frac{1}{2}(1-x_i)^2 \sum_{k\ge 1}   (km-k-1) s_i^k \P(X_i^{(M)}=k).
\end{align*}
 By $n_2\le i \le M-2$, $s_i\ge m-m\delta^3\ge 1$ and   $ \E([(m-1)X_i^{(M)}-1] s_i^{X_i^{(M)}})=0$, we have
 $$\E((((m-1)X_i^{(M)}-1)s_i^{X_i^{(M)}})1_{\{X_i^{(M)}\ge 1\}})\ge \frac{1}{4},$$
  see \eqref{e:-theta+}.    Hence
  $$
  \Delta_i^{(M)}(s)\ge \frac{1}{8}(1-x_i)^2\ge \frac{\delta^2}{128}.
  $$
\qed\\


  Now, we make full prepare for the proof of Proposition \ref{t:pX2}.\\

 {\it Proof of Proposition   \ref{t:pX2}.} Fix $\delta\in (0,\frac{1}{16m})$ and $c_2 :=\E(m^{(1-\delta)X_0^*})<\infty$. Let $i\ge 0$ and $M\ge  i+2$.  Set  $ s_i> 0$  for the value which satisfies
 $ \E([(m-1)X_i^{(M)}-1] s_i^{X_i^{(M)}})=0$ as before.
 we will show that \eqref{e:px2} holds true with
  $c_3 :=128c_2^{m^{n_2+1}} m^{m-1}  \delta^{-2}$, $c_4 :=c_6(m, m\delta^3)$ and $c_5 :=c_7(m, m\delta^3)$.

  Since $X_i$ is stochastically less than $\sum\limits_{j=1}^{m^i} X_{0,j}$ for each $i\ge 0$, we have
$$\E(m^{(1-\delta)X_i})\le \E(m^{(1-\delta)X_0})^{m^i}\le \E(m^{(1-\delta)X_0^*})^{m^i}=c_2^{m^i}.$$
Hence
$$
\prod_{i=0}^{M-1}\E(m^{(1-\delta)X_i})^{m-1}\le \prod_{i=0}^{M-1} c_2^{m^i(m-1)}\le c_2^{m^{M}}.
$$
By definition, $c_2^{m^{n_2+1}}\le c_3$.  So, \eqref{e:px2} is true for $M\le n_2+1$.

If  $M\ge n_2+2$ and if $s_i\le m-m\delta^3$ for some $n_2\le i\le M-2$, then $1\le (s_i\vee 1)\le m-m\delta^3$ and
$\E(((m-1)X_i^{(M)}-1)(s_i\vee 1)^{X_i^{(M)}})\ge \E(((m-1)X_i^{(M)}-1)s_i^{X_i^{(M)}})=0$.
So, we can  apply    Theorem \ref{l:trucation} by replacing $(m, \delta)$ with $(m, m\delta^3)$
to get
$$
\P(X_n=0)\le c_7(m, m\delta^3) e^{-c_6(m, \, m\delta^3)m^{n-M}}, \quad n\ge M.
$$
Such \eqref{e:px2}    holds true in this situation, too.

We are left to prove the case: $M\ge n_2+2$ and  $s_i>m-m\delta^3$ for any $n_2\le i\le M-2$. Let $s\in (1,m]$. Recall the  definition of $\Delta_{i}^{(M)}(s)$ in \eqref{def:Delta}.
Set
$$f_s(k) :=[1-(s-1)k-\frac{(m-1)(m-s)}{m}k(k+1)]s^k,\quad k\in \mathbf{Z}^+.$$
Then $\Delta_{i}^{(M)}(s)=\E(f_s(X_{i}^{(M)}))$. Since $m\ge 2$ and  $s\in(1, m]$,
$$
  -f_s(1)=[-1+(s-1)+2\frac{(m-1)(m-s)}{m}]s=\frac{m-2}{m}(2m-s)s\ge 0.
$$
By observed,  $k\rightarrow -1+(s-1)k+\frac{(m-1)(m-s)}{m}k(k+1)$ and $k\rightarrow s^k$ are  increasing.
Since  $f_s(0)=1$, we also have $-f_s(0)\le -f_s(1)$.  So that,   $-f_s(k)$ is increasing in $k\in \mathbf{Z}^+$.
 Set $X_{i,j}^{(M)} :=X_{i,j}\wedge(M-i)$ for $1\le j\le m$. Since $X_{i+1}=(X_{i,1}+\cdots+X_{i,m}-1)^+$,
 $$
 X_{i+1}^{(M)}=X_{i+1}\wedge (M-i-1)\le (X_{i,1}^{(M)}+\cdots+X_{i,m}^{(M)}-1)^+.
 $$
By the  property of $-f_s(\cdot)$, we obtain $f_s(X_{i+1}^{(M)})\ge f_s((X_{i,1}^{(M)}+\cdots+X_{i,m}^{(M)}-1)^+)$, and so
$$
 \Delta_{i+1}^{(M)}(s)=\E(f_s(X_{i+1}^{(M)})) \ge \E(f_s((X_{i,1}^{(M)}+\cdots+X_{i,m}^{(M)}-1)^+)).
$$

On the other hand, set $H_i^{(M)}(s) :=\E(s^{X_i^{(M)}})$ as before.  Since $X_{i,j}^{(M)}, 1\le j\le m$ are i.i.d copies of $X_i^{(M)}$, we have
$\E(s^{(X_{i,1}^{(M)}+\cdots+X_{i,m}^{(M)}-1)^+})=\frac{1}{s} ~H_i^{(M)}(s)^m+(1-\frac{1}{s})~H_i^{(M)}(0)^m$. So,
 as   (29) in   \cite{collet-eckmann-glaser-martin} and (29) in \cite{bmxyz_questions}, we have
\begin{align*}
&\E(f_s((X_{i,1}^{(M)}+\cdots+X_{i,m}^{(M)}-1)^+))\\
=&\frac{m}{s}\Delta_i^{(M)}(s) H_i^{(M)}(s)^{m-1}-\frac{m-s}{s}[(m-1)sH_i^{(M)'}(s)-H_i^{(M)}(s)]^2 H_i^{(M)}(s)^{m-2}.
\end{align*}
Hence
   \begin{align*}
 \Delta_{i+1}^{(M)}(s) \ge& \frac{m}{s} \Delta_i^{(M)}(s) H_i^{(M)}(s)^{m-1}\\
 &~-\frac{m-s}{s} [(m-1)sH_i^{(M)'}(s)-H_i^{(M)}(s)]^2 H_i^{(M)}(s)^{m-2}.
\end{align*}

 Choose $s=m^{1-\delta}$ now.  Since $\delta\in (0,\frac{1}{16m})$ and $s_i> m-m\delta^3$,
applying Lemma \ref{l:Cauchy_X} gives
$$
[(m-1)sH_i^{(M)'}(s)-H_i^{(M)}(s)]^2\le 2H_i^{(M)}(0)\Delta_i^{(M)}(s)\le 2\Delta_i^{(M)}(s).
$$
Therefore, for $n_2\le i\le M-2$,
\begin{align*}
\Delta_{i+1}^{(M)}(s)
\ge& \frac{m}{s} \Delta_i^{(M)}(s) H_i^{(M)}(s)^{m-1}-2\frac{m-s}{s} \Delta_i^{(M)}(s)   H_i^{(M)}(s)^{m-1}\\
=& \frac{2s-m}{s} \Delta_i^{(M)}(s) H_i^{(M)}(s)^{m-1}.
\end{align*}

By $\delta\in (0,\frac{1}{16m})$ and $m\ge 2$, we have $\frac{2s-m}{s}=2-m^\delta\ge m^{-2\delta}$.  So,
$$
\Delta_{i+1}^{(M)}(s)\ge m^{-2\delta} \Delta_i^{(M)}(s) H_i^{(M)}(s)^{m-1}.
$$
Iterating the above inequalities from $i=n_2$ to $i=M-2$ gives
$$
\Delta_{M-1}^{(M)}(s)\ge m^{-2\delta(M-1-n_2)} \Delta_{n_2}^{(M)}(s) \prod_{i=n_2}^{M-2} H_i^{(M)}(s)^{m-1},
$$
which implies
$$
\prod_{i=n_2}^{M-2} H_i^{(M)}(s)^{m-1}\le  \frac{\Delta_{M-1}^{(M)}(s)}{\Delta_{n_2}^{(M)}(s)} m^{2\delta(M-1-n_2)}.
$$
By  \eqref{e:Delta_low_bound}, we have $\Delta_{n_2}^{(M)}(s)\ge \frac{\delta^2}{128}$. Since $f(x,s)\le f(0,s)=1$ for all $x\in \mathbf{Z}^+$, we have $\Delta_{M-1}^{(M)}(s)=\E(f(X_{M-1}^{(M)},s))\le 1$. Such we obtain
\begin{equation}\label{e:multi_n_2-M}
\prod_{i=n_2}^{M-2} H_i^{(M)}(s)^{m-1}\le \frac{128}{\delta^2} m^{2\delta M}.
\end{equation}
Taking  $\prod_{i=0}^{n_2-1} H_i^{(M)}(s)^{m-1}\le  c_2^{m^{n_2}}$, $H_{M-1}^{(M)}(s)=\E(s^{X_{M-1}\wedge 1})\le s\le m$ and \eqref{e:multi_n_2-M} together, we draw out
\begin{align*}
\prod_{i=0}^{M-1} H_i^{(M)}(s)^{m-1}\le   c_2^{m^{n_2}} \cdot \frac{128}{\delta^2} m^{2\delta M}\cdot m^{m-1}\le c_3 m^{2\delta M}.
\end{align*}
We  have completed the proof of the proposition.
\qed\\

\section{Some further remark}
Recall $P_{X_0}=(1-p)\delta_0+pP_{X_0^*}$ which was defined in \eqref{e:lawX}.
By Remarks \ref{remark1} and \ref{remark2},  under certain situation
$$
\frac{d^k}{dp^k}\lim\limits_{n\rightarrow\infty}\frac{\E(X_n)}{m^n}\Big|_{p=p_c}=\lim\limits_{n\rightarrow\infty}\frac{d^k}{dp^k}\frac{\E(X_n)}{m^n}\Big|_{p=p_c}=0, \quad k\ge 0.
$$
Intuitively, the derivative operation and the limit operation are exchangeable with respect to $\frac{\E(X_n)}{m^n}$ at $p=p_c$.
 It is also interesting to study  $\E(X_n)$ itself, see our previous   papers \cite{xyz_sustainability,xyz_dual}.
We will ask a similar question for $\E(X_n)$ but   under a different definition of $X_0$.

We say a probability measure $\mu$ on $\mathbf{Z}^+$ is subcritical if
$\sum\limits_{k=0}^\infty ((m-1)k-1)m^k\mu(\{k\})<0$, critical if $\sum\limits_{k=0}^\infty ((m-1)k-1)m^k\mu(\{k\})=0$ and supercritical if $\sum\limits_{k=0}^\infty ((m-1)k-1)m^k\mu(\{k\})>0$.
So,  the  choice  of  $P_{X_0}$ which satisfies \eqref{e:lawX} is just
  a linear combination of a subcritical   $\delta_0$ and a supercritical   $P_{X_0^*}$.

   Now we consider another choice of $X_0$.
    Let $\mu$ and $\lambda$ be two  critical probability measures  on $\mathbf{Z}^+$ with $\mu\not=\lambda$. Let $\{X_n, n\ge 0\}$ be the Derrida-Retaux system which satisfies $P_{X_0}=(1-p)\mu+p\lambda$, i.e.,
   $$
  \P(X_0=k)=(1-p)\mu(k)+p\lambda(k), \quad  k\ge 0, \quad p\in (0,1).
  $$
  So that $P_{X_0}$ is a linear combination of two critical probability measures now. It implies $P_{X_0}$ is  critical,  that is to say, $\E(m^{X_0})=(m-1)\E(X_0 m^{X_0})<\infty$ for  $p\in (0,1)$.
    By Theorem A,
  $$
  \lim_{n\rightarrow\infty}\E(X_n)=0.
  $$
 We wonder whether the derivative operation and the limit operation are still exchangeable with respect to $\E(X_n)$, so that, there is a question.
\begin{question} Under  certain integrability condition  for $\mu$ and $\lambda$,
do we have $\lim\limits_{n\rightarrow\infty} \frac{d}{dp}\E(X_n) =0$ for  $p\in (0,1)$ ?
\end{question}

\noindent {\bf Acknowledgments.} We are grateful to Zhan Shi for suggesting the problem and helpful discussion. This work was partially supported by National Natural Science Foundation of China  grant No. 12271351.

\end{document}